\documentclass {article}

\usepackage{PRIMEarxiv}

\usepackage[utf8]{inputenc} 
\usepackage[T1]{fontenc}    
\usepackage{booktabs}       
\usepackage{amsfonts}       
\usepackage{nicefrac}       
\usepackage{microtype}      
\usepackage{lipsum}
\usepackage{fancyhdr}       
\usepackage{graphicx}       
\graphicspath{{media/}}     
\usepackage{amsmath}
\usepackage{soul}
\usepackage{graphicx}
\usepackage[hyphens,spaces,obeyspaces]{url}
\usepackage{xfrac}
\newcommand{\ra}[1]{\renewcommand{\arraystretch}{#1}}
\usepackage{booktabs}
\usepackage{algorithm}
 \usepackage[]{algpseudocode}
\usepackage{tikz}
\usetikzlibrary{arrows,shapes,snakes,automata,backgrounds,petri,fit,positioning, chains}
\usepackage{eqparbox}
\usepackage{url}
\usepackage{hyperref}

\algdef{SE}[SUBALG]{Indent}{EndIndent}{}{\algorithmicend\ }%
\algtext*{Indent}
\algtext*{EndIndent}
 
\usepackage{color}
\usepackage{soul}

\usepackage{subfig}
\usepackage{color, colortbl}
\definecolor{Gray}{gray}{0.9}
\newcommand\correspondingauthor{\thanks{Corresponding author.}}

\usepackage{mwe}
\usepackage{tikz}
\title{A Markov Decision Process Approach for Managing Medical Drone Deliveries
}

\author{
  Amin Asadi \correspondingauthor \\
   Department of Industrial Engineering \\ and Business Information Systems \\
  University of Twente, Enschede, Netherlands \\ \\
    Department of Industrial Engineering\\
  University of Arkansas  \\
  4207 Bell Engineering, Fayetteville\\
  \texttt{amin.asadi@utwente.nl} \\
   \And
  Sarah Nurre Pinkley\\
  Department of Industrial Engineering\\
  University of Arkansas  \\
  4207 Bell Engineering, Fayetteville\\
  \AND
   Martijn Mes \\
   Department of Industrial Engineering \\ and Business Information Systems \\
  University of Twente, Enschede, Netherlands \\
}

\begin{document}
\maketitle

\begin{abstract}
We consider the problem of optimizing the distribution operations at a drone hub that dispatches drones to different geographic locations generating stochastic demands for medical supplies.  Drone delivery is an innovative method that introduces many benefits, such as low-contact delivery, thereby reducing the spread of pandemic and vaccine-preventable diseases.  While we focus on medical supply delivery for this work, drone delivery is suitable for many other items, including food, postal parcels, and e-commerce. In this paper, our goal is to address drone delivery challenges related to the stochastic demands of different geographic locations. We consider different classes of demand related to geographic locations that require different flight ranges, which is directly related to the amount of charge held in a drone battery.  We classify the stochastic demands based on their distance from the drone hub, use a Markov decision process to model the problem, and perform computational tests using realistic data representing a prominent drone delivery company.  We solve the problem using a reinforcement learning method and show its high performance compared with the exact solution found using dynamic programming. Finally, we analyze the results and provide insights for managing the drone hub operations.
\end{abstract}


\section{Introduction} \label{sec: Intro}

During the last decade, there has been substantial growth in the use of drones for various applications, including but not limited to transportation, agriculture, and delivery \cite{DubaiDrone, AgDrones2, UPSDrones, DHL}. Specifically, delivery using drones has received extensive attention as it can reduce air pollution and traffic in congested areas \cite{Dhote2020}. Moreover, drones are a viable option to reach remote locations with inadequate road infrastructure \cite{Malawi19}. During pandemics, drones provide a safe and low contact delivery method, which can effectively slow down the spread of the diseases \cite{Unicef20, McNabb20}. Many companies and organizations, including Vanuatu's Ministry of Health and Civil Aviation \cite{Kent19}, Zipline \cite{Zipline20}, Matternet \cite{Matternet20}, and Manna Aero \cite{ MannaAero20} use drones to deliver and distribute medical supplies such as vaccines, medicine, and blood units. Effective operations of a fleet of drones require addressing battery-oriented issues, including limited flight range,  time-consuming charging operations, and the high price and short lifetime of batteries. \textcolor{black}{In this research, we provide a model that effectively uses the charge inside the batteries to maximize the amount of stochastic demand met. 
Specifically, we classify the stochastic demands according to drones' flight range, which depends on the charge inside drone batteries. We ultimately maximize the expected total met demand for delivering medical items using drones dispatched from a battery swap station located in a drone hub.}

Battery swap stations are a solution to alleviate the aforementioned issues. In battery swap stations, charged batteries are swapped with empty batteries in short minutes. For instance, Matternet provides a station to automatically swap drone batteries used for delivering blood units and medicine between the supplier and hospitals \cite{Matternet20-2}. Besides the quick swapping operation, recharging batteries in anticipation of demand can reduce overcharging and fast charging of batteries, which are shown to accelerate the battery degradation process \cite{Lacey13, Shirk15}. The faster batteries degrade, the quicker the need for battery replacement, leading to a higher cost and environmental waste. The application of a battery swap station is not limited to drones and can be extended to electric vehicles \cite{electrek2018, ChinaDaily19}, electric scooters \cite{Gogoro} and cell phone battery swaps in airports, hotels, and amusement parks \cite{FuelRod}. Notably, the number of electric vehicle swap stations is growing in different regions of the world \cite{ChinaSwap, FranceSwap, IndiaSwap, SlovakiaSwap}. In this research, we consider a \emph{swap station located at a drone hub} that dispatches drones to satisfy multiple classes of stochastic demand for medical supplies, which are classified based on their distance from the station.

Given the growth in the number and applications of battery swap stations, it is crucial to optimally manage the stations' operations to reach the highest performance of the station. Thus, we design a decision-making framework to provide optimal recharging and distribution policies when considering the stochastic demand originating from different geographical locations. 
The drone hub can send drones to locations within their flight range, which differ based on the level of charge inside their batteries. \textcolor{black}{It is a complicated setting, given that the level of charge inside batteries can be any number between 0 and 1, and different combinations of charge levels can be used to satisfy the stochastic demands generated from places located at different distances from the drone hub. Hence, using aggregation and discretization, w}e classify the demand based on the distance between the hub and demand locations. \textcolor{black}{To the best of our knowledge,} we \textcolor{black}{are the first} \textcolor{black}{to use such classification for this problem and link} a level of charge inside batteries with a class of demand such that the demand of each class can be satisfied with batteries having the same or higher level of charge. That is, each class of demand can be satisfied with one or multiple levels of charge inside batteries. We formulate this problem as a stochastic scheduling and allocation problem with multiple classes of demand (SA-MCD).

We use the Markov decision process (MDP) to model the stochastic SA-MCD. It is an appropriate modeling approach for problems like ours that are in the class of sequential decision-making under uncertainty problems \cite{Puterman05}. The decisions are made in discrete points in time or decision epochs. We represent the state of the system in the MDP as the number of batteries within each charge level class.  The actions of the MDP are the number of batteries recharging from one level to an upper level of battery charge. The transition probability is a complex function governed by multiple classes of stochastic demand. In our MDP, the optimal policy determines the maximum expected total reward, which is a function of total weighted met demand of different classes.

 We use backward induction (BI) and a reinforcement learning method with a descending $\epsilon$-greedy exploration feature (D$\epsilon$RL) to solve SA-MCDs. The exploration feature allows the algorithm to take random actions to assist in escaping short-term or local optima by visiting apparently non-attractive states. The descending exploration means that as the algorithm proceeds in iterations, the probability of taking an arbitrary action decreases, and the algorithm perform more exploitation (taking favorable actions) by visiting attractive states. BI can provide exact solutions for problems like SA-MCDs that have finite state and action spaces \cite{Puterman05}. However, BI runs into the curses of dimensionality and faces computational time and memory issues as our problem size increases. Thus, we apply an RL method with an exploration feature that is able to find high-quality approximate solutions for large-scale SA-MCDs, which are not solvable using BI \cite{Powell}. We show the convergence of our RL method and demonstrate its capacity to save computational time and memory. 
 
 
We computationally test the SA-MCD model and solution methods on a case study related to the Zipline drone delivery company, which delivers blood units, vaccines, and medical supplies in Rwanda.  We consider the drone delivery of these supplies from its station located in the Muhanga district, Rwanda, to the reachable hospitals throughout the country.  We consider the population of districts, number of hospitals in each district, number of people using a hospital, and rate of arrivals to each hospital to find the stochastic orders for medical supplies. Then, we convert the orders to the demand for drone missions, given that each drone can carry 2kg of medical products at a time \cite{ZiplineCarry}. For our computational experiments, we classify this demand into two classes based on the distance between the station and each hospital (i.e., closer hospitals are classified in level 1 demand, and farther hospitals are classified in level 2 demand).  We import the real data associated with the distance between locations, the population of districts, flight regulations in Rwanda, and the Zipline drone configuration, including the speed, flight range, and recharging time.

We derive insights from solving SA-MCDs to manage the distribution operations of the swap station using different sets of computational experiments. We provide the optimality gap and average percentage of the met demand using the D$\epsilon$RL method for a modest problem (15-21 drones). Solving the modest size problem shows that the Zipline company needs more drones to satisfy 100\% of the stochastic demand. Hence, we draw the relationship between the number of drones in the station and the amount of met demand using our RL solution for larger instances of SA-MCD. We also analyze the interplay between the different demand classes and the use of higher-level charged batteries to satisfy lower-class demand. 
  
 
 \textbf{\emph{Main Contributions.}}  We summarize the main contribution of this paper as follows. 
 \begin{itemize}
\item We propose stochastic scheduling and allocation problems with multiple classes of demand (SA-MCD) for managing operations of a drone swap station located at drone a hub. We classify the demand based on the distance between the station and hospitals generating the stochastic demands.

\item  We develop an MDP model for the SA-MCD and \textcolor{black}{solve small instances of SA-MCD} using backward induction \textcolor{black}{(BI) and show inability of BI to solve large-scale SA-MCDs required for the real application.}  

\item We \textcolor{black}{propose applying a descending $\epsilon$-greedy reinforcement learning method (D$\epsilon$RL)} to find optimal and near-optimal policies for the station that faces stochastic demand for sending drones to deliver medical supplies. \textcolor{black}{We show high performance of D$\epsilon$RL for solving SA-MCDs.}  

\item We \textcolor{black}{leverage our model and our approximate solution method to be used for a real application, which is} the case study related to Zipline medical supply delivery using drones in Rwanda.

\item We conduct different sets of experiments to derive insights for managing the operations in a swap station to maximize demand satisfaction. We show that demand classification approach in modeling improves demand satisfaction, which is the primary purpose of delivering medical supplies using drones.  
\end{itemize}

The remainder of this paper is organized as follows. In Section \ref{sec: Lit}, we discuss relevant literature related to the modeling, application, and solution methods for the stochastic SA-MCD. In Section \ref{sec: Problem}, we present our Markov Decision Process to model the stochastic SA-MCD. In Section \ref{sec: Sol}, we discuss the exact and approximate solution methods. In Section \ref{sec: Comp}, we outline the computational experiments conducted and provide insights for managing swap station operations. We end with concluding remarks and propose directions for future work in Section \ref{sec: Conc}.

\section{Related work} \label{sec: Lit}

There is a growing interest in the use of drones for various applications. We provide an overview of scientific works, which are more relevant to the model, application, and solution methods presented in this paper. Therefore, we focus on providing an overview of research related to managing operations in swap stations, delivering medical items using drones, Markov Decision Process (MDP) modeling for dynamic problems, demand classification, and reinforcement learning (RL)  methods.  


Many researchers have studied managing swap station operations. Asadi and Nurre Pinkley \cite{Asadi19} present an MDP model to find the optimal/near-optimal policies (number of recharging/discharging and replacement actions) to maximize the expected total profit for the station facing stochastic demands and battery degradation. They solve this problem using a heuristic,  RL methods, and a monotone approximate dynamic programming algorithm \cite{AsadiTs2021} to provide insights for managing the internal operations in the swap stations. Widrick et al. \cite{Widrick16} propose an MDP model for the same problem when no battery degradation is considered. Nurre et al. \cite{Nurre14} do not consider stochasticity and provide a deterministic model to find the optimal policies for managing swap stations. \textcolor{black}{Our work is fundamentally different from the} discussed papers \textcolor{black}{as they} do not consider different demand classes and multiple states of charge of batteries. \textcolor{black}{Besides, our objective is satisfying the amount of met demand, which suits the healthcare delivery application, that differs from the aforementioned papers.}  Kwizera and Nurre \cite{Kwizera2018} propose a two-level integrated inventory model to manage internal operations in a drone swap station delivering to multiple customers (or classes, equivalently) but exclude the uncertainty in the system. To the best of our knowledge, we are the first to introduce the stochastic SA-MCDs for managing internal operations in a swap station facing stochastic demands from different geographical locations.

In recent years, there has been a rapid growth in using drones for many innovative applications \cite{Macrina20}. Several papers \cite{Barmpounakis16, Chang2018, Otto2018, Khoufi2019} review the applications of drones in different contexts, and we refer the reader to Otto et al. \cite{Otto2018} for an extensive review on the optimization approaches for civil applications of drones. Delivering portable medical items such as blood units and vaccines using drones can positively impact the levels of medical service in remote or congested places where roads are not a viable option for transportation and delivery \cite{Otto2018, Dhote2020}.  Several companies are using drones to deliver medical supplies in different parts of the world \cite{Kent19, Zipline20, Matternet20, MannaAero20}.  Notably, we focus on a case study related to Zipline, a drone delivery company that started with 15 drones delivering medical items to remote locations in Rwanda in 2016 \cite{Seeker}.  After successful operations in Rwanda, Zipline expanded its medical delivery service in the south of Ghana using 30 drones in 2019 \cite{Rescue2019}.  During the COVID-19 pandemic, drone delivery provides a fast, cheap, and reliable method to distribute COVID-19 vaccines. For instance, in Ghana, Zipline already delivered 11000 doses and will deliver more than 2.5 million doses in 2021 \cite{Vincent21}.  Draganfly, a Canada-based company, will use drones to distribute COVID-19 vaccines to remote areas of Texas starting in Summer 2021 \cite{Draganfly}.


A drone can store a limited amount of energy that restricts its flight range. This limitation needs to be considered when modeling real-world problems. Common modeling approaches are to use the maximal operation time \cite{Tokekar13, Wang2017} and maximal flying distance \cite{Savuran2016, Guerriero2014}. In this paper, we use the maximal coverage of 80km radius (160km round-trip) \cite{ZiplineFeatures} from our swap station, which is located at the Muhanga Zipline drone hub, for geographically-based demand classification in our case study.

Our SA-MCD problem is in the class of sequential decision-making under uncertainty, and we use a Markov Decision Process (MDP) model, which is appropriate for this class of problems \cite{Puterman05}. There is extensive research on the use of MDPs for stochastic problems in the operations research community. A sample of problems and applications that are close to our research include drone applications \cite{Al-Sabban13, Baek13, Yu15}, dynamic inventory and allocation \cite{Federgruen84, SOMARIN17}, and optimal timing of decisions \cite{millart04, Chhatwal10, Zhang12, Khojandi14}. 

In SA-MCD, we use demand classification that researchers broadly utilize to study scheduling, allocation, supply chain management, and inventory control problems. We discuss a sample of scientific works that used such a classification in combination with a MDP modeling approach. Gayon et al. \cite{Gayon09} provide optimal production policies for a supplier facing multiple classes of demands that are different in the demand rates, expected due dates, cancellation probabilities, and shortage costs. Benjaafar et al. \cite{Benjaafar2011} formulate an MDP model to derive optimal production policies for an assembly system wherein the demands are classified based on the difference between shortage penalties incurred due to the lack of inventory to satisfy orders. Thompson et al. \cite{Thompson2009} categorize patients served by a hospital according to the floors treating patients and the lengths of stay in hospitals. Milnar and Chevalier \cite{Milnar16} use an infinite horizon MDP to model an admission control problem to maximize the expected total profit of a firm serving two classes of customers. The customers are classified based on profit margins, order sizes, and lead time. We use a geographically-based demand classification influenced by different flight ranges of a drone based on the level of charge inside its batteries. An instance of considering travel range for demand classification can be found in \cite{WANG2019} wherein it is assumed that different types of EVs have different drive ranges. They incorporate this information to provide a framework that estimates the charging demands for charging stations and determine the service capacity of the stations without optimizing the system.  As they do not consider optimization, this work is significantly different from our problem that considers optimization under uncertainty.



When an MDP model is large and complex it suffers from the curses of dimensionality  \cite{Powell, Sutton18}, and thus, is traditionally solved using approximate solution methods.  Researchers apply Reinforcement Learning (RL) (or approximate dynamic programming (ADP), the term more used in the operations research community \cite{Powell}) to find approximate solutions that are not solvable using standard exact solution methods, (e.g., dynamic programming \cite{Puterman05}). Examples of various ADP/RL methods in dynamic allocation problems are temporal difference learning \cite{vanr97, Cimen17, Cimen15}, case-based myopic RL \cite{JiangSheng09}, Q-Learning \cite{Chaharsooghi}, value function approximation \cite{Bertmis02, Erdelyi10, Maxwell10}, linear function approximations \cite{Powell05}, and policy iteration \cite{Nasrollahzadeh18}. In this paper, we apply a value function approximation using a look-up table (e.g., \cite{JiangSheng09, Kwon08}) with an $\varepsilon$-greedy exploration feature \cite{Powell, Ryzhov19} to make our RL method visit and update the value of more (both attractive and unattractive) states in the state space. We reduce the exploration rate (increase the exploitation rate) to make the algorithm converge as it proceeds toward iterations. In Section \ref{sec: RL}, we present a comprehensive explanation of our RL method.

\section{Problem description and formulation} \label{sec: Problem}

In this section, we present our Markov Decision Process (MDP) approach to model the stochastic scheduling and allocation problem with multiple classes of demand (SA-MCD).  We proceed by formally describing the classes of demand and the components of our MDP model.  

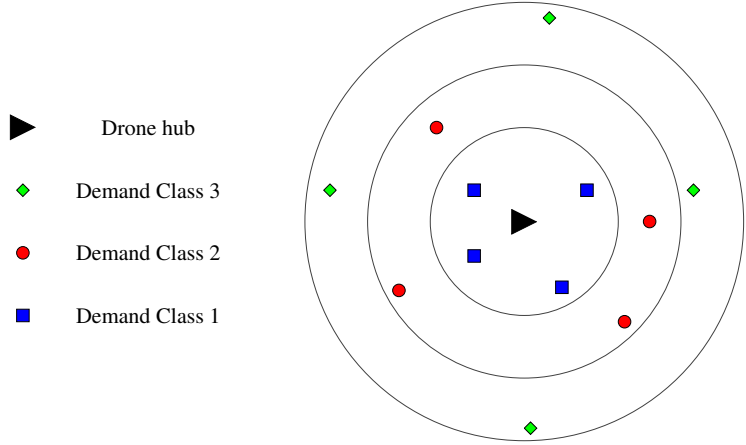
\begin{figure}[ht]
\centering
 \resizebox{0.6\textwidth}{!}{
  \begin{tikzpicture} 
    \draw[darkgray] (0,0) circle[radius=1.5];
    \draw[darkgray] (0,0) circle[radius=2.5];
    \draw[darkgray] (0,0) circle[radius=3.5];

        \draw[black, fill=blue] (0.9,0.4) rectangle(1.1,0.6);
        \draw[black, fill=blue] (-0.9,0.4) rectangle(-0.7,0.6);
        \draw[black, fill=blue] (-0.9, -0.65) rectangle(-0.7,-0.45);
        \draw[black, fill=blue] (0.5, -1.15) rectangle(0.7,-0.95);

        \draw[black, fill=red] (2,0) circle[radius=0.1];
        \draw[black, fill=red] (-1.4,1.5) circle[radius=0.1];
         \draw[black, fill=red] (-2, -1.1) circle[radius=0.1];
         \draw[black, fill=red] (1.6, -1.6) circle[radius=0.1];
           
        \draw[black, fill=green] (2.8,0.5) -- (2.7,0.6) -- (2.6,0.5) -- (2.7,0.4)--(2.8,0.5) ;
        \draw[black, fill=green] (-3,0.5)--(-3.1,0.6) -- (-3.2,0.5) -- (-3.1,0.4)--(-3,0.5);
        \draw[black, fill=green] (0, -3.3)--(0.1, -3.2)--(0.2, -3.3)--(0.1, -3.4)--(0, -3.3);
        \draw[black, fill=green] (0.3, 3.25)--(0.4, 3.35)--(0.5, 3.25)--(0.4, 3.15)--(0.3, 3.25);
        
       \draw[black, fill=black] (-0.2, 0.2)--(-0.2,-0.2)--(0.2,0);

         \node[fill=white] at (-6, 1.5) {Drone hub}; 
       \draw[black, fill=black] (-7.8, 1.5)--(-8.2,1.7)--(-8.2,1.3);   

              
          \node[fill=white] at (-6, -1.5) { Demand Class 1}; 
           \draw[black, fill=blue] (-8.1, -1.6) rectangle(-7.9, -1.4);   
         \node[fill=white] at (-6, -0.5) { Demand Class 2};
          \draw[black, fill=red] (-8, -0.5) circle[radius=0.1];
        \node[fill=white] at (-6, 0.5) { Demand Class 3};
         \draw[black, fill=green] (-8.1, 0.5)--(-8,0.6)--(-7.9, 0.5)--(-8,0.4)--(-8.1, 0.5);
  \end{tikzpicture}
  }
  \caption{An example of demand classification based on the distance between the location of demand and the drone hub.}
  \label{dronehub}
  \end{figure}

For our problem, we consider a set of medical facilities, each with an unknown number of requests (i.e., demand) for drone delivery over time. We know how long a drone needs to fly from the drone hub (located in the swap station) and back to satisfy a request for each medical facility. We cluster medical facilities with similar flight times into demand classes. The demand for each medical facility is then aggregated by demand class. The uncertainty in our MDP is the number of requests (i.e., demand) for each demand class by time. We assume that there is a known probability distribution that governs the uncertainty for each demand class over time. We depict an example of geographically-based demand classification in Figure \ref{dronehub}.

We link each demand class with the required amount of battery charge that is necessary to make the round-trip flight from the drone hub to the medical facility and back. In other words, higher demand classes that are farther from the drone hub require more charge than those closer to the hub. Charging all batteries to full charge ensures that each drone+battery pair can satisfy a request from any demand class.  However, in reality, this strategy results in a higher total cost, longer recharge times, and faster battery degradation. Thus, we make decisions about how many batteries are recharged to different charge levels over time. Our system incorporates time-varying elements, including the mean demand per class over time. We model the system using a finite horizon MDP. 
\textcolor{black}{In addition to the relevance of finite horizon MDP for highly variable and dynamic systems, like ours, the finite horizon model is applicable for the station to check the quality (capacity) of batteries at the end of the time horizon and decide to replace batteries if needed.} We seek to maximize the expected total reward, which equals the summation of the weighted met demand over all demand classes and time periods. We note that the main purpose of our swap station is demand satisfaction; thus, we do not directly incorporate different charging costs and times into our model. However, in the objective function, we use multipliers for the amount of demand of each class that is satisfied using the batteries of higher charge levels.   In this design, higher charge level batteries are capable of satisfying demand for lower classes, but this can be penalized as it caused unnecessary charging costs. \textcolor{black}{Further, because the mission of our drones is delivering medical items, which are often vital at the moment of orders, we do not consider backlogging or rolling over unmet demands.} To maximize the expected total weighted demand, we seek optimal policies indicating how many batteries are charged to each charge level by state and time. We proceed by detailing the specific components of our MDP.

\indent \textit{{Decision Epochs:}} Discrete times within the finite time horizon ${N < \infty}$ in which we make decisions. The set of decision epochs is $T=\{1, 2, \dots, N-1\}.$ 

\indent  \textit{{States:}} The state of the system, $s_t$, is dynamic and defined in the $C$-dimensional state space $S$. Thus, $s_t = (s^1_t, s^2_t, \dots , s^C_t) \in S = (S^1 \times S^2 \times \dots S^C)$, where $S^i$ is the state space for battery charge level $i$ for $i = 1, \ldots, C$. Each $S^i = \{0, \ldots, M\}$ where $M$ represents the total number of batteries.  For each dimension, $s^i_t \in S^i$ equals the number of batteries with $i$ level charge.  The total number of batteries over all charge levels must not exceed $M$ in accordance with Equation \eqref{eq:Si}.  All batteries with a charge level lower than the lowest charge level (i.e., level 1) are implicitly denoted as being level 0 which does not need to be stored but instead, can be calculated as $s^0_t = M - \sum_{i=1}^C s^i_t$.
\begin{equation}
S = \bigg\{(s^1_t, s^2_t, \dots s^C_t): \left(\sum\limits_{i=1}^C s^i_t \le M, \quad \forall t \in T\right)\bigg\}. \label{eq:Si}
\end{equation}

The $C$-dimensional state space corresponds to the $C$ demand classes.  As previously mentioned, each demand class represents a set of medical facilities with similar round-trip delivery flight times.  We link the state space with these demand classes by stating that a drone powered with charge level $i$  is able to satisfy requests from demand class $i$ and lower.  In other words, a request from demand class $i$ is able to be satisfied by a drone powered with charge level $i$ or higher.  We make the assumption that a demand request from class $i$ is satisfied using a battery from the lowest, capable, available charge level.  For example, imagine we have a request from demand class 2.  Any battery with charge level $2, \ldots, C$ is capable of satisfying this request.  If a battery is available with charge level 2, this battery is used.  However, if no batteries are available with charge level 2, then we look to assign a battery with charge level 3, and continue increasing the charge level until a battery is available.  If none are available, we designate this demand as unmet.  With this assumption, we maintain higher charge level battery in inventory.




\indent  \textit{{Actions:}} We use $a_t$ to denote the recharging action at time $t$ using a vector of size ($\sfrac{C(C-1)}{2}$) such that $a^{ik}_t$ represents the number of batteries starting at charge level $i$ which are recharged to level $k$ for $k > i$ at time $t$.  As follows, 
 \begin{equation}\label{eq:at}
a_t = \{a^{ik}_t \in A^{ik}_{s_t}: \  \forall i = 0, 1, \dots, C-1, \ k = i+1, \dots, C \}
\end{equation}
where
\begin{equation}
A^{ik}_{s_t} = \{0, 1, \dots, s^i_t\} \quad \ \ \ \forall i = 1, \dots, C-1, \ k = i+1, \dots, C, \ \text{and} 
\end{equation}
\begin{equation}
A^{0k}_{s_t} = \{0, 1, \dots, M-\sum _{i=1}^{C} s^i_t\}   \ \ \ \forall k = 1, \dots, C. 
\end{equation}

To ensure that the number of batteries selected to be recharged from each charge level does not exceed the number of batteries within that class, we force the actions to satisfy Equations \eqref{eq:limit1} and \eqref{eq:limit2}.
\begin{equation}
\sum_{k=1}^C a^{0k}_t \le M-\sum _{i=1}^{C} s^i_t. \label{eq:limit1}
\end{equation}
\begin{equation}
\sum_{k=i+1}^C a^{ik}_t \le s^i_t  \quad \ \ \ \forall i = 1, \dots, C-1, \ \forall t \in T. \label{eq:limit2}
\end{equation}

In Figure \ref{fig:TransitionExample}, we display the state transitions between different states due to different recharging actions or demand satisfaction for a single battery. Recharging/demand satisfaction increases/decreases the level of charge of a battery depending on the level of recharging/classes of met demand.


\begin{figure}[!ht]
\centering
\begin{tikzpicture}[start chain=going left,node distance=2.4cm,  font=\sffamily\scriptsize]   
\node[state, on chain]                 (C) {$C$};
\node[on chain]                   (g) {$\cdots$};
\node[state, on chain]                 (2) {$2$};
\node[state, on chain]                 (1) {$1$};
\node[state, on chain]                 (0) {$0$};
\draw[
    >=latex,
    auto=right,                      
    loop above/.style={out=75,in=105,loop},
    every loop,
    ]

       (0)   	 edge[loop left] node[left=0.0] {\tiny{Not used}}   (0)
       (1)   	 edge[loop left] node[left=0.0] {\tiny{Not used}}   (1)
       (2)   	 edge[loop left] node[left=0.0] {\tiny{Not used}}   (2)
       (C)   	 edge[loop left] node[left=0.0] {\tiny{Not used}}   (C)

       (0)   	 edge[bend left =50] node[below=0.2] {\tiny{ 1 lvl. charg. from 0 to 1}}   (1)
       (1)   	 edge[bend left =50] node[below=0.2]  {\tiny{1 lvl. charg. from 1 to 2}}   (2)
       (0)  	 edge[bend left =50] node[below=0.1]  {\tiny{2 lvl. charg. from 0 to 2}}   (2)
       (2)       edge[bend left =50] node[below=0.2]  {\tiny{ $C-2$ lvl. charg. from 2 to $C$}}   (C)
       (1)       edge[bend left =50] node[below=0.1]  {\tiny{$C-1$ lvl. charg. from 1 to $C$}}   (C)
       (0)       edge[bend left =50] node[below=0.05]  {\tiny{$C$ lvl. charg. from 0 to $C$}}   (C)

     (g)

       (1)   	 edge[bend left =50] node[above=0.2]  {\tiny{Satisfy Dem. from cls. $1$}}   (0)
       (2)   	 edge[bend left =50] node[above=0.2]  {\tiny{Satisfy Dem. from cls. $1$}}   (1)
       (2)  	 edge[bend left =50] node[above]  {\tiny{Satisfy Dem. from cls. $2$}}   (0)
      (C)   edge[bend left =50] node[above=0.2]  {\tiny{Satisfy Dem. from cls. $C-2$}}   (2)
      (C)   edge[bend left =50] node[above]  {\tiny{Satisfy Dem. from cls. $C-1$}}   (1)
      (C)   edge[bend left =50] node[above]  {\tiny{Satisfy Dem. from cls. $C$}}   (0);

\end{tikzpicture}
\vspace{-0.2in}
\caption{An instance of state transition for a single battery. 
}
	\label{fig:TransitionExample}	
\end{figure}



\indent  \textit{{Transition Probabilities:}} The system transitions from state $s_t$ to a future state $s_{t+1}$ according to the selected action and the realized demand within each demand class.  In our system, the demand at time $t$, denoted $D_t$, is a vector of size $C$, i.e., $D_t = (D^1_t, \dots, D^C_t)$ where each $D^i_t$, for $i=1, \ldots, C$, is a random variable representing the number of requests for demand class $i$.  \textcolor{black}{As we explain later,} the state transitions and the probability transition function are complex, but we illustrate these functions of our MDP where $C=2$.  However, we note that the model could be applied to a problem with $C > 2$.


As our state transition is complex, we first define an intermediate state of the system.  We define the intermediate state of the system as $L= (L^1, L^2)$ and allow this to represent the number of batteries \textcolor{black}{currently charged at levels 1 and 2} after all actions are taken and batteries are allocated to satisfy demand within their class (i.e., batteries with charge levels 1 and 2 are used to satisfy the demand of class 1 and class 2, respectively).  \textcolor{black}{Therefore, the intermediate states should not be mistaken for the post-decision states that show the system's state immediately after making decisions before realizing the uncertainty.} We note, this intermediate transition does not incorporate the batteries from charge level 2 used to satisfy remaining demand from class 1.  The transitions to intermediate  states are governed by Equations \eqref{eq:L1} and \eqref{eq:L2}. In Equation \eqref{eq:L1},  $\text{min} \{s^1_t- a^{12}_t , D^1_t\}$  equals the satisfied demand of class 1 using the available batteries with charge level 1. Similarly, $ \text{min} \{s^2_t , D^2_t\}  $ denotes the satisfied demands of class 2 using batteries with charge level 2 in Equation \eqref{eq:L2}.

\vspace{-0.1in}
\begin{equation}\label{eq:L1}
  \begin{gathered}
L^1= s^1_t +  a^{01}_t - a^{12}_t - \text{min} \{s^1_t- a^{12}_t , D^1_t\}.  
  \end{gathered}
\end{equation}
\begin{equation}\label{eq:L2}
  \begin{gathered}
L^2 = s^2_t +  a^{02}_t + a^{12}_t - \text{min} \{s^2_t , D^2_t\}. 
  \end{gathered}
\end{equation}

Using this intermediate state, we now present the entire state transition equations.  Given $L = (L^1, L^2)$, we can now use the remaining batteries with level 2 charge to satisfy any remaining demand for class 1.  We present the full future state of the system with Equations \eqref{eq:FutureState1} and \eqref{eq:FutureState2}.
\begin{equation}\label{eq:FutureState1}
  \begin{gathered}
s^1_{t+1} = L^1+ \text{min} \big\{\text{max} \{0, D^1_t - (s^1_t- a^{12}_t )\}, \text{max} \{0, s^2_t - D^2_t)\} \big\}.  
  \end{gathered}
\end{equation}
\begin{equation}\label{eq:FutureState2}
  \begin{gathered}
s^2_{t+1} = L^2 - \text{min} \big\{\text{max} \{0, D^1_t - (s^1_t- a^{12}_t )\}, \text{max} \{0, s^2_t - D^2_t)\} \big\}.  
  \end{gathered}
\end{equation}

\begin{figure}[!h]
\begin{center}
\begin{tikzpicture}[shorten >=1pt,draw=black!50, node distance=\layersep, scale=0.5]
\draw (0,0) -- (30,0);
\draw (0.5,2) -- (0.5,-2);
\node [below, label={[align=left, font=\footnotesize]  epoch $t$}] at (2.25,7) {};
\draw [decorate,decoration={brace,amplitude=10pt,raise=90pt},yshift=0pt] (0 ,0) -- (4.5,0) node [black,midway,xshift=0.00cm, yshift= -0.8cm] {};

\node [below, label={[align=left, font=\footnotesize]  Observed  \\ state at \\ epoch $t$}] at (0.5,2) {};
\node [below, label={[align=center, font=\footnotesize]  \# charged \\at level 1 and 2 \\ $(s^1_t, s^2_t)$}] at (1,-5) {};
\node [below, label={[align=center, font=\footnotesize] Decisions for\\ recharging \\ $(a^{01}_t, a^{02}_t, a^{12}_t)$}] at (4.5,2) {};
\draw[->] (4.5,0) -- (4.5,1.5);
\draw (13.25,1) -- (13.25,-1);

\draw [decorate,decoration={brace,amplitude=10pt,mirror,raise=4pt},yshift=0pt] (5,0) -- (13,0) node [black,midway,xshift=0.00cm, yshift=-0.8cm] {};

\draw [decorate,decoration={brace,amplitude=10pt,mirror,raise=4pt},yshift=0pt] (13.75,0) -- (21.5,0) node [black,midway,xshift=0.00cm, yshift=-0.8cm] {};

\node [below, label={[align=center, font=\footnotesize, color=blue] Allocate drones to meet\\ $D^1_t$ and $D^2_t$ \\ using  available \\ batteries with \\ the same level of charge}] at (9.2,-5.5) {};
\node [below, label={[align=center, font=\footnotesize, color=blue] Allocate drones to meet \\ unmet demand of  $D^1_t$ \\ using available \\ batteries with the \\ higher level of charge}] at (18,-5.5) {};
\draw [line width=0.1cm, color=blue] (5,0) -- (12.75,0)  ;
\draw [line width=0.1cm, color=blue] (13.75,0) -- (21.5,0)  ;

\node [below, label={[align=left, font=\footnotesize]  Updated the \\ Intermediate \\ state $(L^1, L^2)$}] at (13.5,2) {};

\draw[->] (22,1.5) -- (22,0);
\node [below, label={[align=left, font=\footnotesize]  epoch $t+1$}] at (26.1,7) {};
\draw [decorate,decoration={brace,amplitude=10pt,raise=90pt},yshift=0pt] (22 ,0) -- (29,0) node [black,midway,xshift=0.00cm, yshift= -0.8cm] {};
\node [below, label={[align=center, font=\footnotesize]   Updated state using \\ info. of charged batteries  \\ and satisfied demand}] at (22,2) {};
\node [below, label={[align=center, font=\footnotesize] Observed  state \\ at epoch $t+1$}] at (28.5,2) {};
\draw (29,2) -- (29,-2);
\node [below, label={[align=center, font=\footnotesize] \# charged \\at level 1 and 2  \\ $(s^1_{t+1}, s^2_{t+1})$}] at (28,-5) {};
\end{tikzpicture}
\caption{Diagram outlining the timing of events for the SA-MCD model.}
\label{fig:eventTiming}
\end{center}
\end{figure}
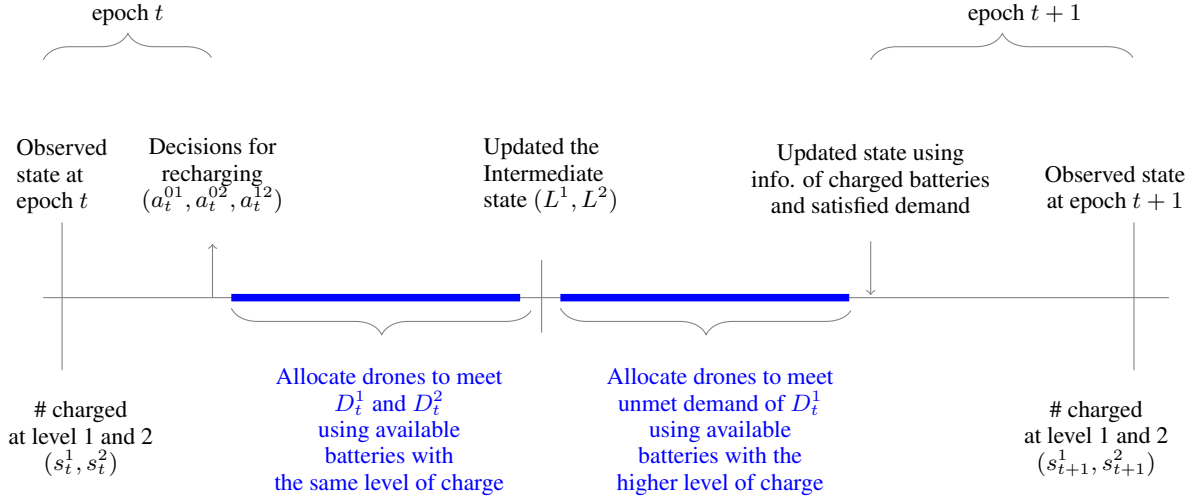

\textcolor{black}{To illustrate the state transitions, we use Figure \ref{fig:eventTiming} to display the timing of events between two consecutive decision epochs. At epoch $t$, we observe the state of system $(s^1_t, s^2_t)$ and make a decision for charging actions $(a^{01}_t, a^{02}_t, a^{12}_t)$. Then, we realize the stochastic demand at epoch $t$. It means the information for the number of demands of each class becomes available after we make the decisions for charging actions. Between two decision epochs, we first allocate drones with level 1 and 2 charges to satisfy the demand class 1 and 2, respectively. Then, we can determine the intermediate state of the system and allocate the leftover batteries of level 2 charge to meet the unmet demand of class 1. Ultimately, at epoch $t+1$, we update the state of system using Equations \eqref{eq:L1}, \eqref{eq:L2}, \eqref{eq:FutureState1}, and \eqref{eq:FutureState2}.}

 In Equations \eqref{eq:FutureState1} and \eqref{eq:FutureState2}, $U^1 = \text{max} \{0, D^1_t - (s^1_t- a^{12}_t )\}$ is the amount of unsatisfied class 1 demand after using level 1 charged batteries and $U^2 = \text{max} \{0, s^2_t - D^2_t)\}$ is the number of leftover level 2 charged batteries after satisfying class 2 demands. Hence, the amount of  remaining class 1 demand that can be satisfied using remaining level 2 charged batteries is the minimum of $(U^1, U^2)$.  We note that our system holds the Markov property that means the system's future state does not depend on the state of the system in the past and can be derived using solely the present state, taken action, and realized uncertainty \cite{Puterman05}. 

 \textcolor{black}{We provide an example to clarify the state transitions. Suppose we have 10 batteries and the system's state is $(s^1_t, s^2_t) = (3,6)$, which means 3 and 6 batteries have level 1 and 2 charges, respectively, and one battery is depleted. 
Lets assume, we take the action $(a^{01}_t, a^{02}_t, a^{12}_t) = (0, 1, 2)$ and the realized demand $(D^1_t, D^2_t) = (5, 2)$. Incorporating this information, the number of met demand of class 1 and 2 using batteries with the same level of charge equals $\text{min} \{s^1_t- a^{12}_t , D^1_t\} = \text{min} \{3-2 , 5\} = 1$ and $\text{min} \{s^2_t , D^2_t\} = \text{min} \{6 , 2\} = 2$, respectively. 
The amount of unsatisfied class 1 demand is $ U^1 = \text{max} \{0, D^1_t - (s^1_t- a^{12}_t )\} = 4$. The number of leftover level 2 charged batteries that can be used to satisfy unmet demand of class 1 is $ U^2 =\text{max} \{0, s^2_t - D^2_t\} = \text{max} \{0, 6 - 2\} = 4$. Using Equations \eqref{eq:L1} and \eqref{eq:L2}, the intermediate state of the system is $(L^1, L^2) = (0, 7)$. Then, we use 4 out of 4 level 2 charged batteries to satisfy the demand of class 1 as $\text{min} (U^1 , U^2) = 4$. These 4 batteries will return to the station with level 1 charge. Hence, the future state $s^1_{t+1} = L^1 +4 = 4$ and $s^2_{t+1} = L^2 -4 = 3$.}  
 

Now, we present the transition probability function from state $s_t$ to state $s_{t+1} = j = (j^1, j^2)$ using Equation \eqref{eq:TRPR}.In this equation, $p^i_x = P(D^i_t = x)$ and $q^i_x = \sum_{u=x}^{\infty} p^i_u = P (D^i_t \ge x) \ \forall i =1, 2$. In all of the cases in Equation \eqref{eq:TRPR}, the intermediate state transitions are calculated using Equations \eqref{eq:L1} and \eqref{eq:L2}. The first case in Equation \eqref{eq:TRPR} calculates the transition probabilities when the stochastic demand of class 1 and 2 is less than the number of charged level 1 and 2 batteries, respectively. The future state equates the intermediate state for each charge level. In the second case, the demand of level 2 is greater than or equal to the number of batteries with level 2 charge; hence the number of batteries with level 2 charge at time $t+1$ equals the number of recently charged batteries from empty or level 1 charge. The future state of the system equates to the intermediate state of the system. The third case is similar to the second case except that in the second case, the stochastic demand of class 1 is less than the number of level 1 charged batteries but in the third case, the demand is greater than or equal to the number of level 1 charged batteries. The fourth case describes the condition that all of the demand for the class 1 charge can be satisfied using all the available level 1 charged batteries plus the leftover batteries of level 2 after satisfying the demand of class 2. The future state of level 2 batteries will be no more than the intermediate state for level 2. The amount of satisfied demand in stage 2 (satisfying demand of class 1 using remaining batteries of class 2) equals the difference between the intermediate and future state of level 2 charged batteries. The fifth case is similar to the fourth case, except that all of the demands of class 1 can not be satisfied in the first and second stage of the demand satisfaction process.

\begin{equation}\label{eq:TRPR}
p(j | s_t,  a_t)=
    \begin{cases}
(p^1_ {s^1_t +  a^{01}_t - a^{12}_t - j^1})(p^2_{s^2_t +  a^{02}_t + a^{12}_t  - j^2}) & \parbox[t]{0.5\textwidth}{ \text{if } $a^{01}_t < L^1 \le s^1_t +  a^{01}_t - a^{12}_t,\\  a^{01}_t + a^{12}_t  < L^2 \le s^2_t +  a^{02}_t + a^{12}_t,\\
 j^2 = L^2\text{, and } j^1 = L^1$} \\  \\
(p^1_ {s^1_t +  a^{01}_t - a^{12}_t - j^1})(q^2_{s^2_t})& \parbox[t]{0.5\textwidth}{\text{if } $a^{01}_t < L^1 \le s^1_t +  a^{01}_t - a^{12}_t,\\ a^{01}_t + a^{12}_t  = L^2, j^2 = L^2\text{, and } j^1 = L^1$}  \\ \\
(q^1_ {s^1_t  - a^{12}_t })(q^2_{s^2_t}) & \parbox[t]{0.5\textwidth}{\text{if } $a^{01}_t = L^1 ,  a^{01}_t + a^{12}_t  = L^2,\\ j^2 = L^2\text{, and } j^1 = L^1$} \\ \\
(p^1_ {s^1_t - a^{12}_t + L^2 - j^2})(p^2_{s^2_t +  a^{02}_t + a^{12}_t  - L^2}) & \parbox[t]{0.5\textwidth}{\text{if } $a^{01}_t = L^1, a^{01}_t + a^{12}_t  < L^2 \le s^2_t +  a^{02}_t + a^{12}_t,\\ a^{02}_t + a^{12}_t  < j^2 \le L^2\text{, and } j^1 = a^{01}_t + L^{2}_t  -j^2$} \\ \\
(q^1_ {s^1_t - a^{12}_t + L^2 - j^2})(p^2_{s^2_t +  a^{02}_t + a^{12}_t  - L^2}) & \parbox[t]{0.5\textwidth}{\text{if } $a^{01}_t = L^1, a^{01}_t + a^{12}_t  < L^2 \le s^2_t +  a^{02}_t + a^{12}_t ,\\ a^{02}_t + a^{12}_t  = j^2\text{, and } j^1 = a^{01}_t + L^{2}_t  -j^2$}  \\ \\
0 & \text{otherwise.} 
    \end{cases}
\end{equation}

We note that if the number of classes increases to $C=3$, then there are 13 different cases to be considered. Further, we need to consider two intermediate states to capture demand satisfaction with the same level of charge and one and two-level of charge higher. For the application of our problem, two demand classes are adequate to observe a significant improvement in the quality of solutions (see Section \ref{sec: Comp}) without making the model excessively complicated. If the application requires having more than three classes, we suggest using alternative strategies such as non-Markovian modeling and/or only applying ADP/RL/Q-learning methods, which does not require an explicit transition probability function.

\textcolor{black}{Despite the apparent complexity of our transition probability function, we would like to point out the benefit of our MDP model.  Incorporating the assumptions for the order of demand satisfaction and timing of events in our model enables us to capture the system's dynamic without including the past state(s) information. As a result, our model holds the Markov property, allowing us to benefit from the MDP properties and guaranteed solution methods.}

\indent  \textit{{Reward:}} We calculate the immediate reward of taking action $a_t$ when the transition from state $s_t$ to state $s_{t+1} = j$ occurs using Equation \eqref{eq:immidiateReward}
\begin{equation}\label{eq:immidiateReward}
\begin{gathered}
r_t (s_t, a_t, j) = \rho^{11}(s^1_t + a^{01}_t - a^{12}_t - L^1) +   \rho^{21} (L^2 - j^2) +  \rho^{22} (s^2_t + a^{02}_t + a^{12}_t - L^2 )
\end{gathered}
\end{equation}

\noindent for $ t = 1, \dots, N-1$. We use $\rho^{ij}$ to put weights on the amount of met demand of class $j$ using level $i$ charged batteries. 

In the first term, $(s^1_t + a^{01}_t - a^{12}_t - L^1)$ denotes the number of  level 1 charged drones used to satisfy class 1 demand. In the second term, $(L^2 - j^2)$ equals to the number of level 2 charged batteries used to satisfy class 1 demand. In the third term, $(s^2_t + a^{02}_t + a^{12}_t - L^2 )$ determines the number of level 2 charged batteries used to satisfy class 2 demand. We note that our objective is to maximize the expected total satisfied demand which does not directly incorporate cost. However, with adjusting the weights of $\rho^{ij}$, we implicitly include a cost factor via assigning a penalty/reward to demand satisfaction with an excessive level of charge. For instance, when $\rho^{21} = \rho^{22} = \rho^{11} = 1$, there is no benefit in recharging batteries to level 1 to satisfy class 1 demand because level 2 charged batteries can satisfy class 1 demand by generating the same reward while these batteries can satisfy class 2 demands, too. However, if $\rho^{21} = 0.5$, then we expect to recharge to/use more level 1 charged batteries to satisfy class 1 demand given that $\rho^{11} = 2 \rho^{21}$, which means we penalize the reward of satisfying class 1 demand with level 2 charged batteries.  We vary this parameter and analyze the results in Section \ref{sec: Comp}. In practice, the drone delivery companies can select the value of this parameter based on the insights or their preferences.

At the end of the time horizon we calculate the terminal reward. We assume that no action is taken at the end of the time horizon and that all remaining batteries can be used to satisfy future demand, and there is sufficient demand for each level.  Thus, we define the terminal reward using Equation \eqref{eq:terminalReward}.
 \begin{equation}\label{eq:terminalReward}
r_N (s_N) =
 \rho^{11} s^1_N + \rho^{22} s^2_N  
\end{equation}

We calculate the immediate expected reward $r_t (s_t, a_t)$, using the immediate reward and transition probability functions given by Equation \eqref{eq:expectedReward},
\begin{equation}\label{eq:expectedReward}
\begin{gathered}
r_t (s_t, a_t) =  \sum_{j, L \in S} \bigg[ p_t(j | s_t, a_t) \big( \rho^{11}(s^1_t + a^{01}_t - a^{12}_t - L^1) +   \rho^{21} (L^2 - j^2) +  \\ \rho^{22} (s^2_t + a^{02}_t + a^{12}_t - L^2 ) \big)  \bigg]. 
\end{gathered}
\end{equation} 
We derive the decision rules, $d_t(s_t) : s_t \rightarrow A_{s_t}$ , from the action set to maximize the total expected reward. Because we select a single action based on the present state, which does not depend on the past states and actions, our decision rules belong to the Markovian decision rules \cite{Puterman05}. A policy $\pi$ is a sequence of decision rules for all decision epochs, that is $d^{\pi}_t (s_t) \ \ \forall \ t \in T$. We can calculate the expected total reward of policy $\pi$ for the problems starting from an arbitrary initial state $s_1$ using Equation \eqref{eq:totalReward}. The optimal policy, $\pi^*$, maximizes the expected total reward.  
\begin{equation}\label{eq:totalReward}
\begin{gathered}
V^{\pi}_N (s_1) =  \mathbb{E}^{\pi}_{s_1} \bigg[\sum_{t=1} ^ {N-1} r_t (s_t,a_t) + r_N (s_N)\bigg].
\end{gathered}
\end{equation}

\section{Solution methodology} \label{sec: Sol}
In this section, we first present the exact solution method, backward induction (BI), and continue with our approximate solution method, the reinforcement learning method with a descending $\epsilon$-greedy exploration feature (D$\epsilon$RL) to solve the stochastic SA-MCDs.

\subsection{Backward induction} \label{sec: BI}

As our Markov Decision Process (MDP) model has finite state and action spaces, there is at least one deterministic optimal policy \cite{Puterman05}; thus, backward induction (BI) can determine such policy (number of recharging actions) that maximizes the expected total reward or weighted met demand over time. Let $V^*_t (s_t)$ be the optimal value function equivalent to the maximum expected total reward from decision epoch $t$ onward when the system is in state $s_t$. Then, we can use the optimality (Bellman) equations, given by Equation \eqref{eq: u}, to find the optimal policies for all the decision epochs when moving backward in time. That is,  BI sets the value of being in state $s_N$ at the end of the time horizon $N$ to be equal to the terminal reward value given by Equation \eqref{eq:terminalReward}.  Then, the algorithm starts from the last decision epoch and finds the optimal actions ($a^*_{s_t,t}$) and corresponding values ($V^*_t(s_t)$) using Equations \eqref{eq: u} and \eqref{eq: Arg} stepping backward in time. The algorithm aims to find the optimal expected total reward over the time horizon, $V_1^*(s_1)$, for state $s_1$, which is the system's initial state at time $t=1$. In other words, solving the optimality equations for $t=1$ is equivalent to the expected total reward over the time horizon. 
\begin{equation}
V^*_t(s_t)=\max_{a_t\in A_{s_t}}\left\{r_t(s_t,a_t)+\sum_{j\in S}p_t(j \mid s_t, a_t)V_{t+1}(j)\right\}.  
\label{eq: u}
\end{equation}
\begin{equation}
a^*_{s_t,t} = \text{arg max}_{a_t\in A_{s_t}}\left\{r_t(s_t,a_t)+\sum_{j\in S}p_t(j|s_t, a_t)V_{t+1}(j)\right\}.
\label{eq: Arg}
\end{equation}

The size of state space, action space, transition probability, and optimal policies (the best actions for all the states over time) are functions of $O(M^2)$ and $O(M^3)$, $O(M^7N)$, $O(M^2N)$, respectively. As the size of the problem increases, it becomes challenging for BI to find the optimal solution due to the curses of dimensionality, which causes a drastic increase in computational time and memory. Hence, we proceed with presenting our reinforcement learning (RL) method, which is capable of circumventing such problems  \cite{Powell, Sutton18}.

\subsection{Reinforcement learning} \label{sec: RL}

In this section, we explain the reinforcement learning method to find near-optimal policies and to overcome the curses of dimensionality \cite{Powell} of the stochastic scheduling and allocation problem with multiple classes of demand (SA-MCD). We use a value iteration approach,  reinforcement learning with a descending $\epsilon$-greedy exploration feature (D$\epsilon$RL), which provides high-quality approximate solutions for SA-MCD ( see Section \ref{sec: Comp}). We proceed by introducing the notation and continue with our RL features and procedure. 

\begin{table}[!htbp]
\ra{1.3}
\footnotesize
\centering
\scalebox{1}{
\begin{tabular}{@{}ll@{}}
\hline
Notation &					Description \\ \hline
$\tau_1$ & 				The number of core RL iterations\\
$\tau_2$ & 				The number of sample paths of demands (realized uncertainty)\\
${V}^n_t (s_t)$ & 			The optimal value of being in state $s_t$ at time $t$ for iteration $n$ \\
$\overline{V}^n_t (s_t)$ & 		The approximate value of being in state $s_t$ at time $t$ for iteration $n$ \\
$\hat{\upsilon}_t^n(s_t)$ & 	The observed value of state $s_t$ at time $t$ for iteration $n$ \\
$\alpha_n$ & 				The step-size value at iteration $n$\\ 
$\varepsilon^n$ & 				The exploration rate at iteration $n$ \\ 
$u$ & A random number that is used to select exploration or exploitation\\
$z_t^n(s_t^n)$ & 			The smoothed value of being in state $s_t$ at time $t$ for iteration $n$ \\ \hline
\end{tabular}}
\caption{Notation used in the reinforcement learning algorithm.}
\label{tab:Notation}
\end{table}

We first determine the number of drones ($M$), decision epochs ($N-1$), ($\tau_1$) iterations, and $(\tau_2)$ sample paths in our RL method. Then, we initialize the approximate value at the end of the time horizon, $t = N$, for all iterations using the terminal reward function given by Equation \eqref{eq:terminalReward}. For every iteration, we select an initial state $s^n_1$. To select the action, we use the $\varepsilon$-greedy method \cite{Powell} that allows exploring the action space works as follows. We generate a random number, $u$. Then, we compare $u$ with the exploration rate, $\varepsilon^n$, at iteration $n$. We use a descending function for the exploration rate over the iterations to ensure more states (both attractive and unattractive) are visited and facilitate the algorithm convergence. If $Rand < \varepsilon^n$, we select a feasible action arbitrarily. Otherwise, we generate $\tau_2$ sample paths of demands (realized uncertainty) and select the action that maximizes the observed value $\hat{\upsilon}_t^n(s_t)$ according to Equations \eqref{eq:Arg3} and \eqref{eq:uOberveration}, wherein $\overline{V}^{n-1}_{t+1} (s_{t+1})$ is used to approximate the value of $\mathbb{E}({V}_{t+1} (s_{t+1}) \mid s_t, a_t)$ for each sample path. If an action, is selected over multiple sample paths, we use the average of $\overline{V}^{n-1}_{t+1} (s_{t+1})$ as the approximation. The observed value and the approximated value at the previous iteration are smoothed using a stepsize function. This value is now used as the present approximation value of the observed state. When an action is selected, we sample an observation of uncertainty (generate a realized value for stochastic demand) to find the future state. The algorithm steps forward in time and moves to the future observed state until it reaches the last decision epoch and new iteration starts. The same process is repeated until $\tau_1$ iterations are completed.  
\begin{equation}
 a^n_{s_t,t} =
\underset{a_t\in A_{s_t}}{\arg\max}\left\{r_t(s_t,a_t) + \overline{V}^{n-1}_{t+1} (s_{t+1})\right\}. 
\label{eq:Arg3}
\end{equation}
 
\begin{equation}
\hat{\upsilon}_t^n(s_t^n)=
  \max_{a_t\in A_{s_t}}\left\{r_t(s_t,a_t) + \mathbb{E}({V}_{t+1} (s_{t+1}) \mid s_t, a_t)\right\}. 
\label{eq:uOberveration}
\end{equation}

\begin{algorithm}[ht]
\caption{Reinforcement Learning Method with a Descending $\epsilon$-greedy Exploration Feature (D$\epsilon$RL)}
\label{alg: RL1}
\begin{algorithmic}[1]
  \footnotesize
\State Initialize $M$ drones, $N-1$ decision epochs, $\tau_1$ iterations, and $\tau_2$ sample paths
\State Set $\overline{V}^n_N(s) = r_N(s)$ for $s \in S$ and $n = 1, \ldots, \tau_1$  \State Set $n = 1$
\While{$n \leq \tau_1$}
\State Select initial state $s_1^n$
	\For{$t = 1, \ldots, N-1$}
	\State Generate a random number $u$
	\If {$Rand < \varepsilon^n$} 
		\State Sample an observation of the uncertainty, $D_t$
		 \State Determine a random feasible action, $a^t_n$ 
		 \State Calculate the observed value, $\hat{\upsilon}_t^n(s_t^n)$
	\Else 
	\State Generate $\tau_2$ sample paths of the uncertainty
	\State Select an action that maximizes $\hat{\upsilon}_t^n(s_t)$ over $\tau_2$ sample paths
	\State Sample an observation of the uncertainty, $D_t$
	\State Calculate the observed value, $\hat{\upsilon}_t^n(s_t^n)$
	\EndIf	
		\State Smooth the new observation with the previous approximated value, \begin{equation} z_t^n(s_t^n) = ({1 - \alpha_n})\overline{V}^{n-1}_t(s_t^n) + {\alpha_n}\hat{\upsilon}_t^n(s_t^n) \label{eq:Smoothing} \end{equation}
				\State Update the present approximation using the smoothed value, 
				$\overline{V}^{n}_t(s_t^n)  \leftarrow z_t^n(s_t^n) $

		\State Determine next state, $s_{t+1}^n$
	\EndFor
	\State Increment $n = n+1$
\EndWhile
\end{algorithmic}
\end{algorithm}


\section{Computational results}  \label{sec: Comp}

In this section, we explain the results of solving the stochastic scheduling and allocation problems with multiple classes of demand (SA-MCD) using realistic data related to the drone delivery company Zipline.  We created this dataset to mimic the geographical locations of the Zipline drone hub and hospitals in Rwanda, Africa, the population of districts, flight regulations in the country, Zipline drone configuration, including the speed, flight range, and recharging time. We solve modest SA-MCDs (15-21 drones) using exact solution methods.  \textcolor{black}{We note that we use the number of batteries and the number of drones interchangeably as each drone has a battery pack, which might be consist of multiple batteries that can be charged in parallel.} As we run into the curses of dimensionality for larger instances of SA-MCD, we present the results of our reinforcement learning (RL) method that can provide near-optimal solutions for the modest instances and solve larger problem instances. We deduce managerial insights for managing the swap station's distribution operations that maximize the expected total weighted met demand of multiple classes. We proceed by first explaining the data.  


\subsection{Data}\label{data}

In this paper, we present an actual case study related to the vital operations of Zipline drones, which are delivering medical items in Rwanda, Africa. The Zipline station is located at the Muhanga district, west of Rwanda's capital city, Kigali. We focus on drone delivery to satisfy the stochastic demand for blood units originating from hospitals in Rwanda.  
In Table \ref{DataSet}, we summarize the input data, including, the name of each hospital, their location, the distance between the station and each hospital, the approximated population of people using each hospital, the number of blood units and flights needed per day for each hospital, and the demand class associated with each hospital. We categorize the demand into two classes based on the distance between the Zipline station and hospitals. As the flight range of the drone is estimated to be 80km \cite{Zipline80km}, we assume demands of hospitals located within [0km, 40km) and [40km, 80km] fall into class 1 and class 2, respectively. The demand class \emph{NA} means that the hospital is not reachable from the Zipline station and excluded from further analysis. We exclude six hospitals (denoted by \emph{NA} in the last column of Table \ref{DataSet}) from the 33 identified hospitals as they are located at a distance out of the drone's flight range from the origin, Zipline station. We consider ten hospitals in class 1 and 17 hospitals in class 2 which are located throughout 15 distinct districts in Rwanda. 

\begin{table}[!ht]\centering
\ra{0.8}
\scalebox{0.65}{
\begin{tabular}{@{}llrrrrcc@{}}\toprule
Hospital 	&	District	&	Dist. to Zipline	&	Pop. reach 	&	Pop. need 	&	Pop. need 	&	 Rounded \# of flights	&	Class of	\\
name	&		&	 station (km)	&	the hospital	&	blood unit/year	&	blood unit /day	&	needed/day	&	demand	\\ \midrule
Nyamata  	&	Bugesera 	&	34.3	&	361914	&	7238.3	&	19.8	&	10	&	1	\\
Butaro  	&	Burera 	&	73.5	&	336582	&	6731.6	&	18.4	&	10	&	2	\\
Nemba  	&	Gakenke 	&	47.8	&	169117	&	3382.3	&	9.3	&	5	&	2	\\
Ruli  	&	Gakenke 	&	27.6	&	169117	&	3382.3	&	9.3	&	5	&	1	\\
Kiziguro 	&	Gatsibo 	&	75.8	&	216510	&	4330.2	&	11.9	&	6	&	2	\\
Ngarama 	&	Gatsibo 	&	77.5	&	216510	&	4330.2	&	11.9	&	6	&	2	\\
Byumba 	&	Gicumbi 	&	61.3	&	395606	&	7912.1	&	21.7	&	11	&	2	\\
Gakoma	&	Gisagara 	&	34.9	&	161253	&	3225.1	&	8.8	&	5	&	1	\\
Remera Rukoma 	&	Kamonyi 	&	20.3	&	340501	&	6810.0	&	18.7	&	10	&	1	\\
Kibuye Referral 	&	Karongi 	&	48.2	&	110603	&	2212.1	&	6.1	&	4	&	2	\\
Kirinda  	&	Karongi 	&	22.6	&	110603	&	2212.1	&	6.1	&	4	&	1	\\
Mugonero  	&	Karongi 	&	55.6	&	110603	&	2212.1	&	6.1	&	4	&	2	\\
Gahini  	&	Kayonza 	&	85.4	&	172079	&	3441.6	&	9.4	&	5	&	NA	\\
Rwinkwavu  	&	Kayonza 	&	94.0	&	172079	&	3441.6	&	9.4	&	5	&	NA	\\
Kirehe  	&	Kirehe 	&	99.6	&	340368	&	6807.4	&	18.7	&	10	&	NA	\\
Kabgayi 	&	Muhanga	&	5.2	&	319141	&	6382.8	&	17.5	&	9	&	1	\\
Kibungo  	&	Ngoma 	&	85.1	&	336928	&	6738.6	&	18.5	&	10	&	NA	\\
Muhororo &	Ngororero 	&	22.6	&	333713	&	6674.3	&	18.3	&	10	&	1	\\
Shyira &	Nyabihu 	&	46.0	&	294740	&	5894.8	&	16.2	&	9	&	2	\\
Nyagatare 	&	Nyagatare 	&	104.9	&	465855	&	9317.1	&	25.5	&	13	&	NA	\\
Kaduha 	&	Nyamagabe 	&	40.8	&	170746	&	3414.9	&	9.4	&	5	&	2	\\
Kigeme 	&	Nyamagabe 	&	53.8	&	170746	&	3414.9	&	9.4	&	5	&	2	\\
Kibogora 	&	Nyamasheke 	&	77.5	&	190902	&	3818.0	&	10.5	&	6	&	2	\\
Nyanza 	&	Nyanza 	&	31.9	&	323719	&	6474.4	&	17.7	&	9	&	1	\\
Munini 	&	Nyaruguru 	&	76.8	&	294334	&	5886.7	&	16.1	&	9	&	2	\\
Kabaya 	&	Rubavu 	&	44.7	&	201831	&	4036.6	&	11.1	&	6	&	2	\\
Gitwe 	&	Ruhango 	&	22.4	&	159943	&	3198.9	&	8.8	&	5	&	1	\\
Ruhango 	&	Ruhango 	&	19.8	&	159943	&	3198.9	&	8.8	&	5	&	1	\\
Kinihira 	&	Rulindo 	&	50.7	&	143841	&	2876.8	&	7.9	&	4	&	2	\\
Rutongo	&	Rulindo 	&	41.9	&	143841	&	2876.8	&	7.9	&	4	&	2	\\
Mibilizi 	&	Rusizi 	&	107.3	&	200429	&	4008.6	&	11.0	&	6	&	NA	\\
Murunda &	Rutsiro 	&	48.3	&	324654	&	6493.1	&	17.8	&	9	&	2	\\
Rwamagana 	&	Rwamagana 	&	73.7	&	313461	&	6269.2	&	17.2	&	9	&	2	\\ \bottomrule
\end{tabular}
}
\caption{The data associated with blood unit delivery using Zipline drones in Rwanda, Africa.}
\label{DataSet}
\end{table}

\begin{figure}[!ht]
\begin{center}
\includegraphics[scale=0.35]{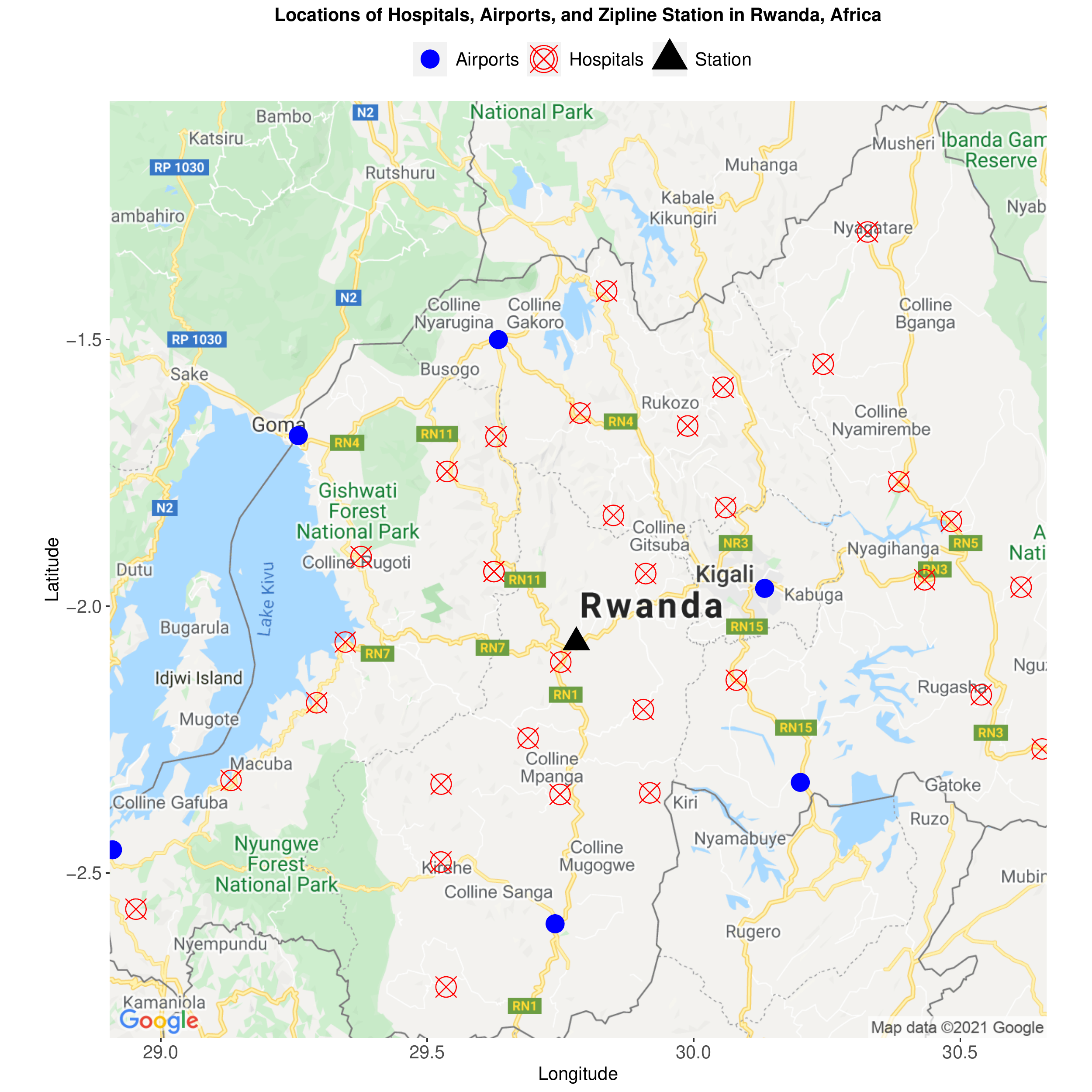}
\caption{Locations of hospitals (demand nodes), the swap station located in the Zipline drone hub, and airports in Rwanda.}
\label{Map}
\end{center}
\end{figure}

We approximate the air travel distance between the station and hospitals using the Haversine formula \cite{Sinnott84} that is broadly used to find the distance between two points on the earth. When calculating the distance, we consider the rules for flying drones in Rwanda, which does not allow drones to fly within a 10km radius from airports \cite{DroneAirport}. Therefore, we need to adjust the travel distance between the station and two hospitals, Kiziguro and Rwamagana. Hence, we calculate the closest travel distances such that the flights to these destinations do not violate the rules for flying drones in Rwanda. In Figure \ref{Map}, we display the geographical locations of airports, hospitals, and the Zipline station. 

Consistent with \cite{Swartzman70, Armony2015}, we use a non-homogenous Poisson process to determine the patients' arrival to hospitals. We examine the daily operations of the drone swap station wherein the time between two consecutive decision epoch is 90 minutes (1.5hr). Thus, $N = \sfrac{24}{1.5} + 1 = 17$. The 90-minute intervals provide adequate time for drones to receive charge \cite{DroneRecharge} and complete a round-trip from the furthest delivery mission to the station given that the maximum speed of the drone is 127km/hr  \cite{DroneSpeed}.

\begin{figure}[!ht]
\begin{center}
\includegraphics[scale=0.22]{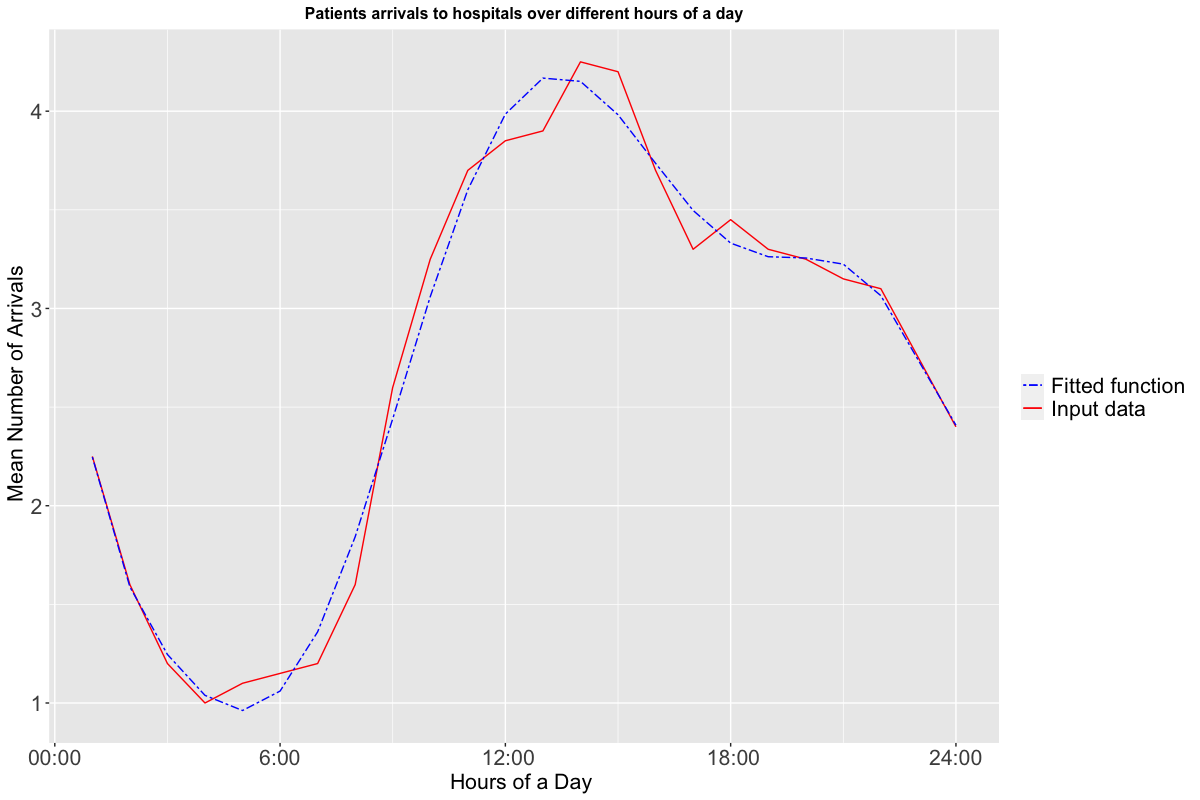}
\caption{Pattern of patients arrivals to a hospital over a day.}
\label{FittedFunction}
\end{center}
\end{figure}

To derive the mean demand of blood units per time $t$, we apply the following process. First, we determine the number of people using a particular hospital based on the population of each district. If more than one hospital is located in a district, we evenly distribute the district's total population over the number of hospitals located in that district. Second, we calculate the number of blood units needed per year for each hospital by multiplying an estimated portion of the population that needs blood units per year (2\% recommended by the World Health Organization (WHO) \cite{BloodEstimate}) by the number of people using that hospital.  The yielded number is an overestimate of the number reported by \cite{BloodOverEstimate}. However, we use this number to account for pessimistic cases wherein the station faces more demand. Third, we divide the number of blood units required per year by 365 to find the number of blood units needed per day for each hospital. Next, we use the pattern of patient arrivals to hospitals, consistent with \cite{Green07, Tiwari2014, Jones2007}, to derive the mean demand for blood units of time $t$ over a day. The pattern in the literature indicates an ascending trend of arrivals from 6:00am to the peak at noon, followed by a descending trend from noon to 6:00am of the following day. Specifically, we use the data from \cite{Green07} and fit a polynomial function for generating the mean arrival rate of time $t$. In Figure \ref{FittedFunction}, we display the mean demand of \cite{Green07} and our fitted function. We scale the mean demand of blood units of time $t$ such that the summation of the scaled demand over a day equals the calculated number of blood units required per day for each hospital.  Then, as each drone can carry two units of blood \cite{ZiplineCarry}, we divide the mean demand of blood units of time $t$ by two to find the mean demand for flights for each hospital. Finally, the mean demand for either demand class, $\lambda^1_t, \lambda^2_t$, is the aggregation of mean demands for flights from the hospitals within each demand class.  

In the first experiment, we consider Zipline has a fleet of 15 drones \cite{Seeker}.  We set $\rho^{11} = \rho^{22} = 1$ and $\rho^{21} = 0.5$ indicating that satisfying a demand of class 1 using a level 1 charged battery (partially-charged) generates more immediate reward than satisfying that demand using a level 2 charged battery (fully-charged). This setting implies the company provides less reward when drones with excessive level of charge are used to satisfy demands, which can be interpreted as a penalty to account for unnecessary higher recharging costs incurred.  

The setting of our  D$\epsilon$RL parameters is as follows. The number of core iterations is $\tau_1 = 200000$. As we will see later, the algorithm will converge after 50000 iterations; however, as the computational time is in a matter of minutes, we keep 200000 iterations for our experiments. We test $\tau_2 = 1, 5, 10, \dots, 50$ and observe that increasing the parameter to the value of 30 reduces the optimality gap and increases the robustness of the results and computational time, but excessive increase in the value of $\tau_2$ only magnifies the computational time with little improvement in the quality of the result; thus, we set $\tau_2 = 30$. We use the adaptive stepsize function provided by George and Powell \cite{George06}. We use $\varepsilon^n = \sfrac{1}{n}$ to adjust the value of the exploration rate at iteration $n$ used in the $\varepsilon$-greedy approach to select policies within our RL method. With this function, we ensure a higher rate of exploration/exploitation in early/late iterations, which is desirable for visiting more states and enabling the algorithm to converge as it proceeds with each iteration.

\subsection{Discussion and analysis}\label{DA}

In this section, we first feed the data explained in Section \ref{data} to solve the problem using exact and approximate solution methods, Backward Induction (BI) and the reinforcement learning  method with a descending $\epsilon$-greedy exploration feature (D$\epsilon$RL), respectively. Then, we analyze the optimal policies (BI solutions) and assess the quality of near-optimal solutions derived from D$\epsilon$RL. Moreover, as the drone delivery company can control and adjust $\rho^{21}$ (the weight of satisfying class 1 demand using level 2 charged batteries), we analyze the impact of changing the parameter's value on the station's operations and amount of met demand.  We also conduct different sets of experiments to solve instances of the problem and answer the following questions.  How many batteries are needed in the station to satisfy a certain level of the stochastic demand? What is the contribution of classifying the demand on the demand satisfaction? We use a high-performance computer with four shared memory quad Xeon octa-core 2.4 GHz E5-4640 processors and 768GB of memory for running all of our computational tests. 

\subsubsection{Comparing results of BI and RL}

In this section, we present the results from solving stochastic SA-MCD with BI and D$\epsilon$RL and using the data presented in Section \ref{data}. The system's initial state is $s_1 = (0,15)$, which means all 15 batteries are charged to level 2 (fully-charged). The time horizon is one day and $N=17$ wherein the first decision epoch is at midnight and the time between any two consecutive decision epochs is 90 minutes. That is, the decisions are made at 16 decision epochs, $t$, where $t$ = 00:00, 1:30, 3:00, \dots, 10:30, 12:00, 13:30, \dots, 22:30. The results of D$\epsilon$RL in the last iteration can differ in terms of the converged value and met demand. Hence, we generate independent sample paths to report robust results and evaluate our RL method's performance. Our preliminary experiments reveal that our result is robust when $>$100 sample paths of demand are incorporated for evaluation (the percentage change in the converged values is satisfactory and less than 0.1). Hence, to yield even a higher robustness level, we calculate the average percentage of met demand using Equations \eqref{ch3: DemOverTime} and \eqref{ch3: AvgDemOver} wherein we generate 500 sample paths of realized demand using the Poisson distribution at time $t$ for each class of demand. 
\begin{multline}
\text{(\%) of Demand Met over Time for a Sample Path} = \\ \bigg|\frac{ \text{Tot. \#   Met Dem. over Time - Tot. \#   Realized Dem. over Time}}{\text{Tot. \#  Realized Dem. over Time}}\bigg| * 100\%.
\label{ch3: DemOverTime}
\end{multline}
\begin{multline}
\text{Average (\%) of Demand Met} = \\ \frac{ \sum\limits_{i=1}^{\text{\#  of Sample Paths}} \text{(\%) of  Met Demand  over Time for Sample Path $i$ }}{\text{\#  of Sample Paths}}.
\label{ch3: AvgDemOver}
\end{multline}

\begin{figure}[ht]
\begin{center}
\includegraphics[scale=0.35]{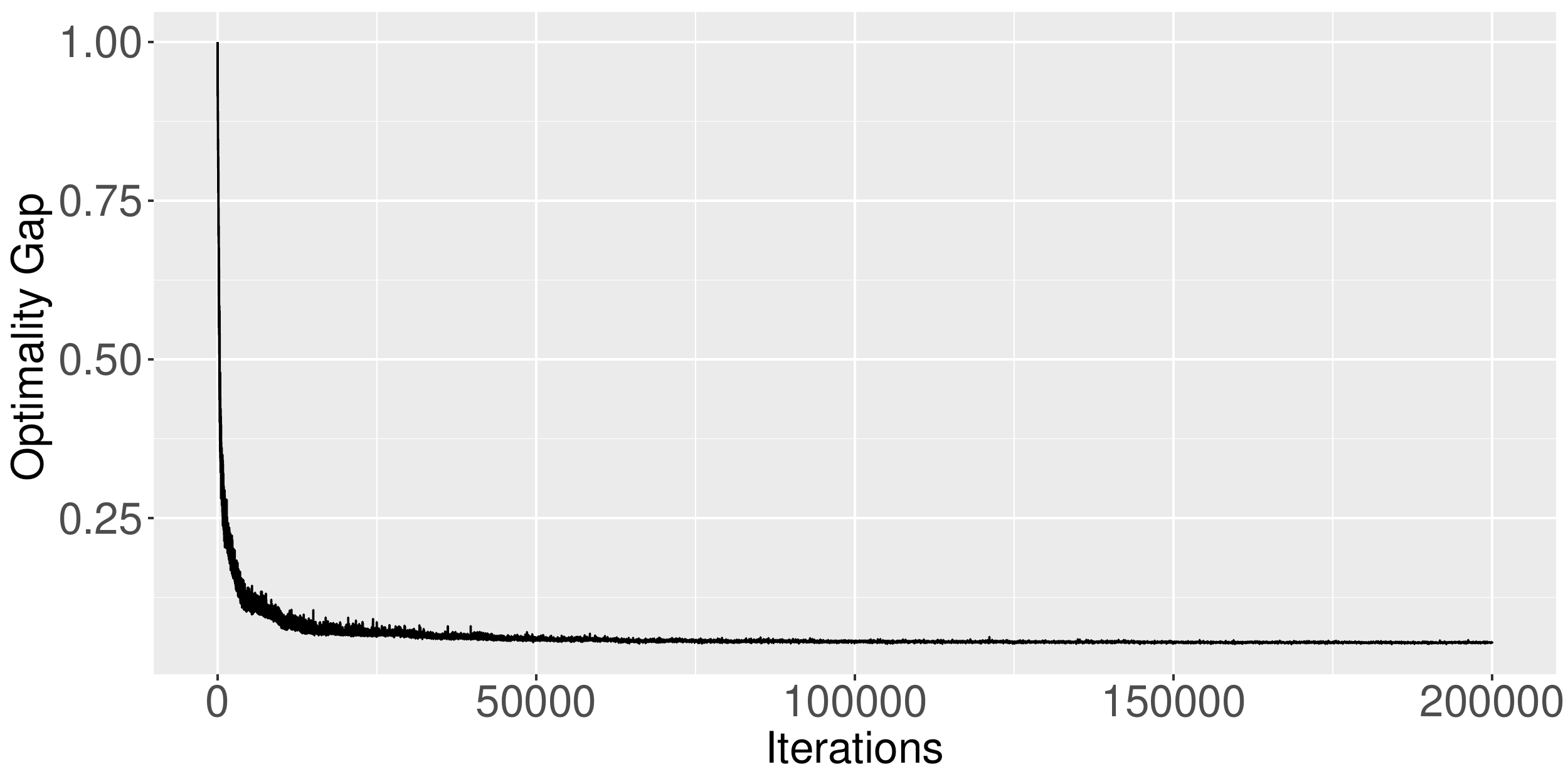}
\caption{Expected total reward convergence of D$\epsilon$RL.}
\label{Convergence15}
\end{center}
\end{figure}

We summarize the results in Table \ref{BasicBI-RL}. The percentage of met demand over a sample path equals the total number of demand met over the total realized demand of both classes. We report the average of 500 sample paths in Table \ref{BasicBI-RL}. We provide more detailed results about demand satisfaction by class later in Table \ref{Rhos}. As shown, D$\epsilon$RL is faster and can generate a high-quality solution with 5.3\% of optimality gap (derived from Equation \eqref{ch3: OptGap}) in 8 minutes. In Figure \ref{Convergence15}, we show the convergence of our proposed D$\epsilon$RL method. 

\begin{table}[H]\centering
\ra{1.2}
\scalebox{0.7}{
\begin{tabular}{@{}ccccc@{}}\toprule
Solution & &Expected  & Average Met  &  Computational \\
Method & &Total Reward & Demand  (\%)  & Time (s) \\ \cline{1-1} \cline{3-5}
BI & & 115.1  & 63.7  & 6740.7\\
D$\epsilon$RL & & 109.0 & 60.9 & 483.7\\
Benchmark & & 105.6 & 58.5 & 6442.1 \\
 \bottomrule
\end{tabular}
}
\caption{The expected total reward and computational time of solving SA-MCD with $\rho^{21} = 0.5$ and $M=15$ using BI and RL.}
\label{BasicBI-RL}
\end{table}

\begin{equation}
\text{Optimality  Gap} = \bigg|\frac{\text{Exp. Tot. Reward BI - Exp.Tot. Reward D$\epsilon$RL Method} }{\text{Exp. Tot. Reward BI}}\bigg| * 100\%.
\label{ch3: OptGap}
\end{equation}

\textcolor{black}{In Table \ref{BasicBI-RL}, we also provide the results of an easy-to-implement intuitive benchmark. As the model's objective is to maximize the amount of met demand, our proposed benchmark policy makes all the empty batteries fully-charged (usable for demand satisfaction of all classes) at every decision epoch. We note the benchmark can be quickly implemented, but it needs a considerable computational time to calculate the exact expected total reward. \textcolor{black}{The computational time for evaluating the benchmark includes all of the operations of backward induction, except the \textit{for loop} for finding the best actions as the actions are derived from the (fixed) benchmark policy.} Our results show that D$\epsilon$RL outperforms the benchmark policy regarding the expected total reward, average demand met, and computational time. Thus, we only continue with D$\epsilon$RL as the superior approximate solution method.}

In Table \ref{tab:BIRL}, we provide the summary of the results from solving the problem with 15 to 21 drones using BI and D$\epsilon$RL. We note, BI cannot find the optimal solution when the number of drones is greater than 21. For 15-21 drones, D$\epsilon$RL provides an optimality gap of less than 6\% for all instances and significantly reduces computational time. The maximum difference between the average percentage of met demand of D$\epsilon$RL and BI is less than 5\%. The results indicate the high performance of our D$\epsilon$RL method in providing approximate solutions.

\begin{table}[!ht]\centering
\ra{0.8}
\scalebox{0.7}{
\begin{tabular}{@{}cccrrccccc @{}}\toprule
	& & \multicolumn{3}{c}{  Backward Induction (BI)} & & \multicolumn{4}{c}{  Reinforcement Learning (D$\epsilon$RL)} \\   \cline{3-5} \cline{7-10} \addlinespace[0.2cm]
$M$	& & 	Average Met  & Computational & Memory		 & & Average Met  & Computational & Memory  & Optimality\\
	& & 	Demand (\%) & Time (s) 		& Used (GB) 	& & Demand (\%)   & Time (s)  & Used (GB) & Gap (\%) \\ \midrule
	
15	& & 63.7 & 6740.7  &  34.2		&  	& 60.9 & 483.7  & 13.4  & 5.3  \\ 
16	& & 66.9 & 11630.4 &  68.1		&  	& 64.3 & 494.1  & 13.4  & 3.3 \\ 
17	& & 69.9 & 19356.4 &  125.8		&  	& 65.2 & 517.0  & 13.4  & 5.0 \\ 
18	& & 72.6 & 31405.8 &  192.2		&  	& 70.6 & 534.5  & 13.4  & 3.4 \\ 
19	& & 75.4 & 49812.5 &  287.5		&  	& 72.6 & 575.7  & 13.4  & 3.5 \\ 
20	& & 77.9 & 77170.6 &   411.9		&  	& 73.6 & 598.1  & 13.4  & 4.8 \\ 
21	& & 80.2 & 117154.0 &   621.1		&  	& 78.3 & 600.1  & 13.4  & 2.7 \\ 
 \bottomrule
\end{tabular}
}
\caption{Computational time, memory used, and average percentage of met demand over time for 500 sample paths when $\rho^{21} = 0.5$ using BI and D$\epsilon$RL methods.}
\label{tab:BIRL}
\end{table}

We can find the visited states and policies over time (sample paths of visited states and policies) using the sample paths of realized demands. Given the initial state, taken actions, realized demand, and the state transition functions given by Equations \eqref{eq:FutureState1} and \eqref{eq:FutureState2}, we can find the future visited state. The consecutive visited states form the sample path of states. The sample paths of policies are the consecutive selected actions derived from the BI and D$\epsilon$RL solution methods.

\begin{figure}[ht]
\begin{center}
\includegraphics[scale=0.4]{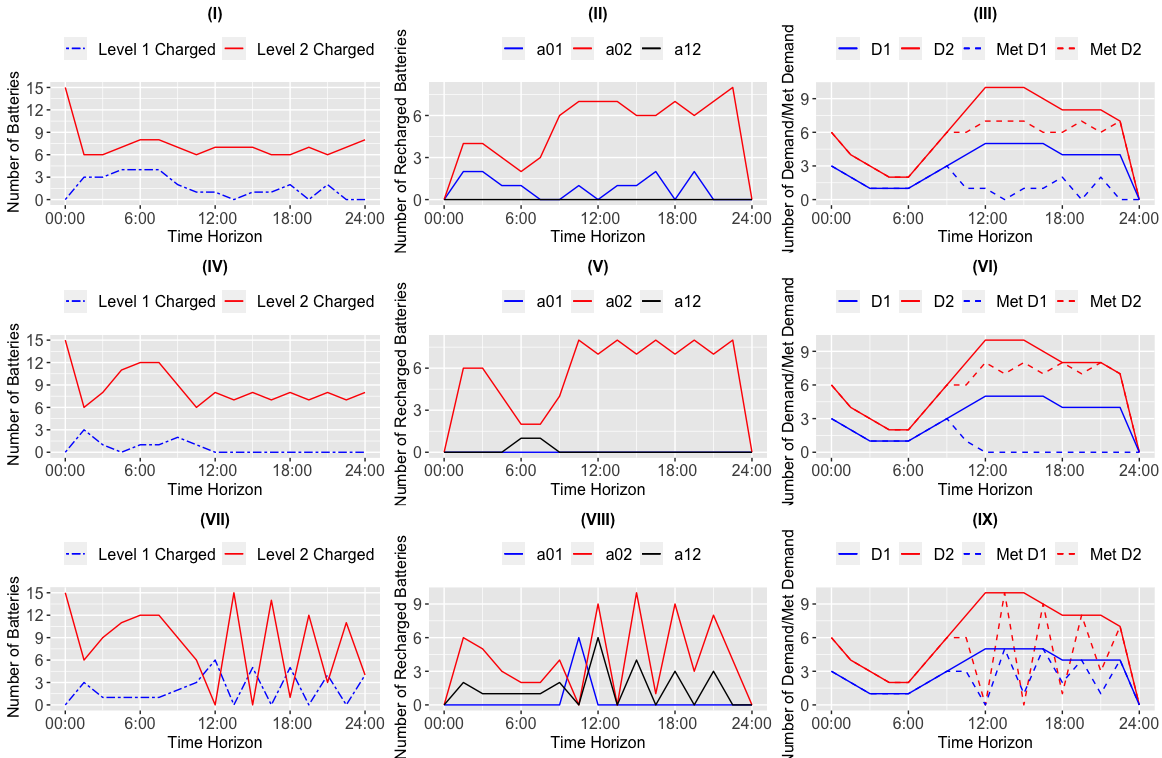}
\caption{Sample paths of states (I, IV, VII), optimal policies (II, V, VIII), and demands (III, VI, IX) for $\rho^{21} = 0.5, 1, 2$ when the realized demand of either class equal mean demand.}
\label{Fig4May2021}
\end{center}
\end{figure}

In Figure \ref{Fig4May2021}, the first, second, and third row depict sample paths of states, optimal policy, and met demands when the stochastic demand equals mean demand at time $t$ and $\rho^{21} = 0.5, 1, \text{and} \ 2$, respectively, for $M=15$. Intuitively, when $\rho^{21} = 0.5$, the level 1 charged batteries are used to satisfy class 1 demand and $a^{12}_t = 0$. When $\rho^{21}$ increases more batteries of class 1 are recharged to class 2 to be used for satisfying the demand of class 1. When $\rho^{21} = 1$, no batteries are recharged to level 1 because there is no difference between the value of satisfying the class 1 demand using level 1 and level 2 charged batteries. Hence, the optimal policy includes recharging batteries to level 2 to be used for either classes of demand. For all values of $\rho$, we observe more recharging actions when demands are in peak periods.


We also compare the average amount of met demand of either class and the optimal policies over time when $\rho^{21}$ (controlled parameter by the station) varies between 0.5 to 2 with 0.1 increments. We calculate the percentage of demand satisfaction of either class by inputting the associated value of each class into Equations \eqref{ch3: DemOverTime} and \eqref{ch3: AvgDemOver} for different values of $\rho^{21}$. We use Equation \eqref{ch3: Action} to find the average number of actions over time over all sample paths.  We summarized the result based on 500 sample paths of realized demand in Table \ref{Rhos}. \textcolor{black}{We note that we display the \textit{number} of actions in the last three columns of Table \ref{Rhos}, but we use \textit{percentage} for the other columns related to demand satisfaction.} Intuitively, as $\rho^{21}$ increases, more/less class 1 demand is satisfied using drones with level 2/level 1 charged batteries. On average, more recharging occurs from level 1 to level 2 to satisfy the demand of either class. For smaller values of $\rho^{21}$, increasing the parameter value provides more incentive for recharging more drones up to level 2 and satisfying both demand levels. However, for larger values of $\rho^{21}$, as level 2 charged batteries are used more to satisfy class 1 demand, fewer level 2 charged drones are available to satisfy class 2 demand. 
\begin{multline}
\text{Average Number of Action $a^{01}$,  $a^{02}$, and $a^{12}$} = \\ \frac{\text{Total Number of the Action over Time over Sample Paths}}{\text{Total Number of Sample Paths}}  * 100\%.
\label{ch3: Action}
\end{multline}

\begin{table}[h]\centering
\ra{0.8}
\scalebox{0.7}{
\begin{tabular}{@{}ccccccccccc @{}}\toprule
$\rho^{21}$&	& Avg Met   	         &  Avg Met  		&	Avg Met 			& Avg Met  		& Avg Met 	&  	&	Avg   & Avg   & Avg \\
		  &	& Demand  C1 with    & Demand  C1 with & 	Demand  C1(\%)	& Demand  C2(\%)  	& Both Class (\%) 	& 	&	$a_{01}$  & $a_{02}$  & $a_{12}$ \\
		  &	& L1 Charge  (\%)      & L2 Charge  (\%)   	&					& 			  	& 				&  	&	 		&   		   &  		 \\ \cline{1-1} \cline{3-7} \cline{9-11} \addlinespace[0.2cm]
0.5		  & 	& 30.0  				& 13.8 		&	43.6				& 		74.5		& 63.6 			& 	& 0.67 		& 	4.97    & 0.05 \\
0.6		  & 	& 25.7  				& 16.2 		&	41.9				& 		76.2		& 64.2 			& 	& 0.45 		& 	5.14    & 0.06 \\
0.7		  & 	& 23.8  				& 17.6 		&	41.4				& 		76.8		& 64.4 			& 	& 0.35 		& 	5.21    & 0.08 \\
0.8		  & 	& 22.3  				& 18.9 		&	41.3				& 		77.1		& 64.5 			& 	& 0.27 		& 	5.25    & 0.09 \\
0.9		  & 	& 21.0  				& 20.0 		&	41.1				& 		77.1		& 64.5 			& 	& 0.20 		& 	5.29    & 0.11 \\
1.0		  & 	& 15.9  				& 24.9 		&	40.9				& 		77.5		& 64.6 			& 	& 0.00 		& 	5.32    & 0.27 \\
1.1		  & 	& 13.9  				& 27.0 		&	41.9				& 		77.5		& 64.6 			& 	& 0.00 		& 	5.24    & 0.41 \\
1.2		  & 	& 13.1  				& 27.6 		&	40.7				& 		77.5		& 64.5 			& 	& 0.00 		& 	5.22    & 0.46 \\
1.3		  & 	& 12.6  				& 28.0 		&	40.7				& 		77.5		& 64.5 			& 	& 0.00 		& 	5.20    & 0.49 \\
1.4		  & 	& 12.2  				& 28.5 		&	40.7				& 		77.3		& 64.4 			& 	& 0.00 		& 	5.17    & 0.53 \\
1.5		  & 	& 11.5  				& 29.7 		&	41.3				& 		76.6		& 64.1 			& 	& 0.02 		& 	5.08    & 0.62 \\
1.6		  & 	& 12.4  				& 33.4 		&	45.8				& 		73.4		& 63.7 			& 	& 0.15 		& 	4.76    & 0.84 \\
1.7		  & 	& 12.5  				& 36.4 		&	48.9				& 		70.4		& 62.9 			& 	& 0.25 		& 	4.48    & 1.04 \\
1.8		  & 	& 11.1  				& 38.1 		&	49.2				& 		69.7		& 62.5 			& 	& 0.28 		& 	4.34    & 1.17 \\
1.9		  & 	& 10.3  				& 39.2 		&	49.5				& 		69.0		& 62.1 			& 	& 0.31 		& 	4.24    & 1.26 \\
2.0		  & 	& 9.1 				& 40.4  		&	49.5				& 		68.6		& 61.9 			&  	& 0.32 		&      4.15    & 1.35 \\
 \bottomrule
\end{tabular}
}
\caption{Average met demand and policies over time for 500 sample paths for different values of $\rho^{21}$.}
\label{Rhos}
\end{table}

A significant finding is that 15 drones are not sufficient to satisfy the demand of either class. We proceed with analyzing the impact of increasing the number of drones in the station on the amount of met demand. 

\subsubsection{Analysis on the number of required drones }

In this section, we solve the problem for a larger number of drones in the station to find the relationship between this number and the amount of met demand. The analysis provides significant insights for drone delivery companies given the high price to purchase and maintain drones in swap stations. First, we note that backward induction (BI) can solve the problem with at most 21 drones using our computational resources. We summarized the amount of met demand, computational time, and memory used to solve the problem for 15 to 21 drones in Table \ref{tab:BIRL} using BI and RL. 
For $M > 21$, we report the results of our D$\epsilon$RL method in Table \ref{tab:RL}.

\begin{table}[!ht]\centering
\ra{0.8}
\scalebox{0.7}{
\begin{tabular}{@{}ccrrcccrrc @{}}\cline{1-5} \cline{7-10} \addlinespace[0.2cm]
$M$	& & 	Average Met  & Computational & Memory		 & & $M$ & Average Met  & Computational & Memory \\
	& & 	Demand (\%) & Time (s) 		& Used (GB) 	 & &	   & Demand (\%)   & Time (s)  & Used (GB)  \\  \cline{1-5} \cline{7-10} \addlinespace[0.2cm]
21	&  & 78.3 & 599.2 & 13.4  & & 41 &   96.6 & 1703.0 &13.4 \\ 
22	& &  78.9 & 654.1 & 13.4  & & 42 &   97.4 & 1761.6 &13.4  \\ 
23	& &  81.0 & 690.1 & 13.4  & & 43  & 97.8 & 1819.7 &13.4  \\  
24	& &  82.0 & 705.2& 13.4  &  &  44 &   98.3 & 1835.9 &13.4 \\  
25	& &  84.5 & 765.1 &13.4  &   &  45 &   98.6 & 1969.3 &13.4 \\ 
26	& & 86.7 & 781.4 & 13.4  &   & 46 &   98.7 & 2106.8 &13.4  \\  
27	& &   87.3 & 935.2 & 13.4  &   &  47  &   99.6 & 2227.0 &13.4 \\   
28	& &   88.5 & 1036.8 & 13.4  &   &  48 &   99.6 & 2582.8 &13.4 \\  
29	&  &   91.2 & 1037.8 & 13.4  &   &  49 &   99.7 & 2717.7 &13.4 \\   
30	& &   91.7 & 1221.1& 13.4  &   & 50 &   99.7 & 2976.0 &13.4 \\  
31	& &   91.8 & 1298.3 & 13.4  &   &  51 &   99.9 &2984.6 &13.4  \\ 
32	&  &   94.2 & 1330.8& 13.4  &   &  52 &   99.9 & 3166.9 &13.4  \\   
33	&  &   94.2 & 1337.3 & 13.4  &   &  53 &   99.9 & 3169.9 &13.4  \\   
34	&  &   94.9 & 1374.1 & 13.4  &   &  54 &  100.0 & 3289.2 &13.4  \\  
35	& &   94.9 & 1384.5 & 13.4  &   &  55 &  100.0 & 3314.5 &13.4   \\ 
36	& &   95.6 & 1406.9 & 13.4  &  &   56  &  100.0 & 3379.8 &13.4  \\   
37	& &   95.7 & 1426.6  & 13.4  &   &  57 &  100.0 & 3431.3 &13.4  \\  
38	& &   95.8 & 1453.5 & 13.4  &   &  58 &  100.0 & 3489.0 &13.4   \\  
39 	& & 	95.9 & 1633.4 & 13.4  &    & 59 & 100.0  & 3593.0 &13.4   \\  
40	& &   96.1 & 1651.4 & 13.4  &   & 60 & 100.0 & 3597.5 & 13.4  \\ 

 \cline{1-5} \cline{7-10}
\end{tabular}
}
\caption{Computational time, memory used, and average percentage of met demand over time for 500 sample paths when $\rho^{21} = 0.5$ using the D$\epsilon$RL method.}
\label{tab:RL}
\end{table}

\begin{figure}[H]
\begin{center}
\includegraphics[scale=0.26]{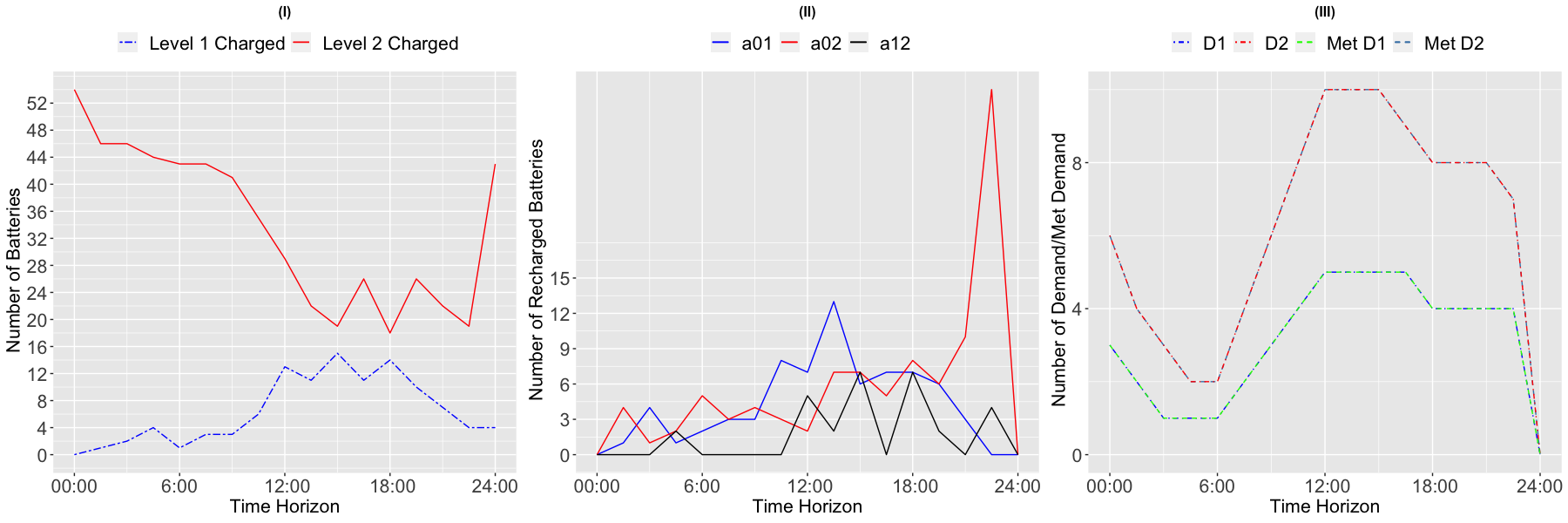}
\caption{Sample paths of states, optimal policies, and demands for 54 drones, $\rho^{21} = 0.5$ when the realized demand of either class equal mean demand.}
\label{Poicy54}
\end{center}
\end{figure}

As shown, when $M \ge 54$, the average percentage of met demand over time for 500 sample paths is 100\%. We depict a sample path of policies for $M=54$ when demand equals mean demand in Figure \ref{Poicy54}. As all batteries are initially available with a level 2 charge, we recharge fewer drones in the early morning. Then, we recharge more batteries from 6:00 to 18:00 as the demand of either class increases. Overall, more batteries are recharged to level 2 to satisfy the demand of either class. 


\subsubsection{Demand Classification Contribution}

In this section, we compare the outputs of the models with and without demand classification to illustrate the contribution of classifying the stochastic demand. We focus on demand satisfaction as the crucial metric to assess the station's success in delivering medical supplies. We provide this metric for a different number of drones, which is important for decision-makers given the high price of purchasing and maintaining drones. 

\begin{figure}[h]
\begin{center}
\includegraphics[scale=0.4]{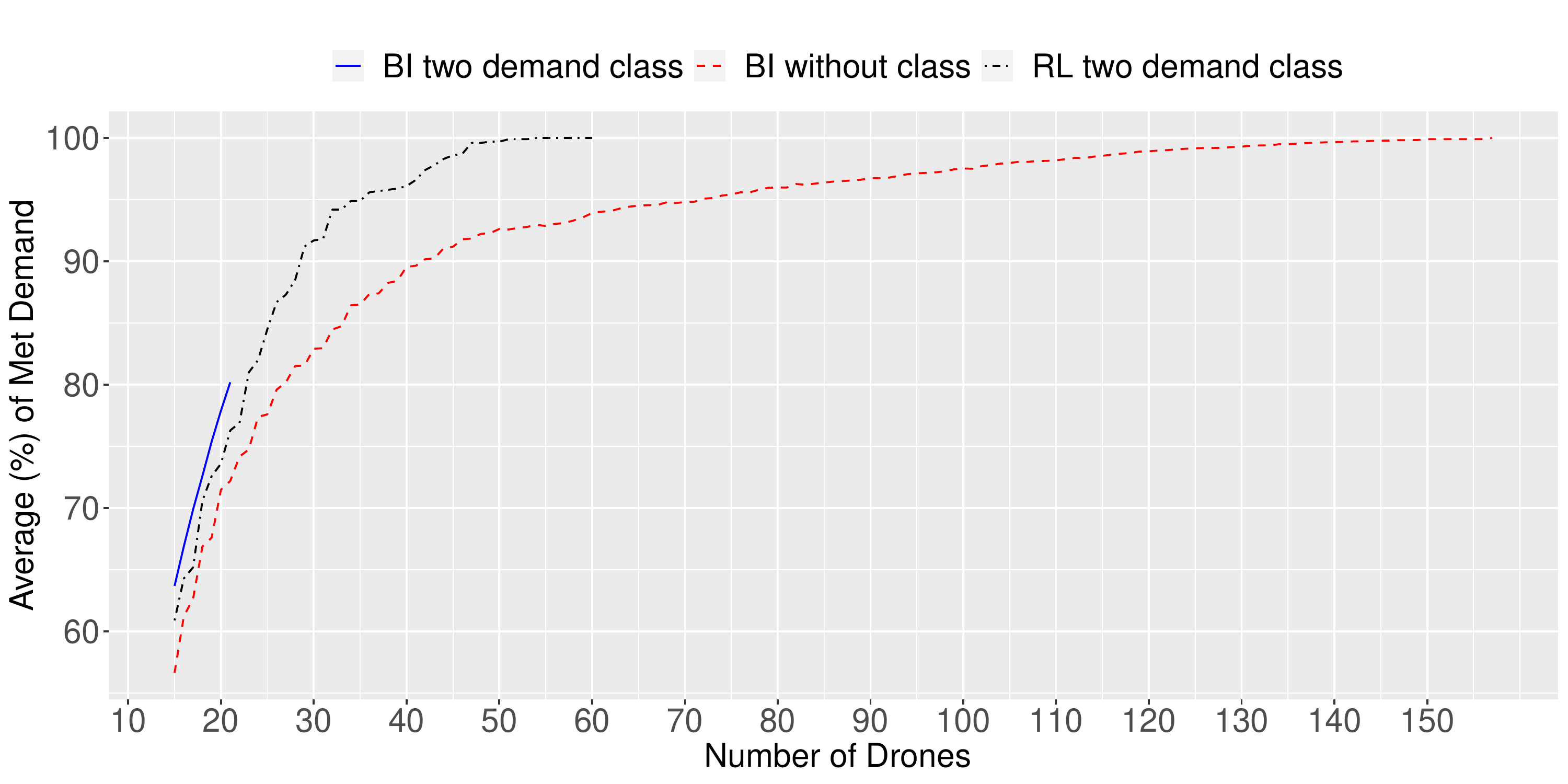}
\caption{Average percentage of met demand for 500 sample paths when $\rho^{21} = 0.5$ using different models and solution methods.}
\label{ClassVsNoClass}
\end{center}
\end{figure}

In Figure \ref{ClassVsNoClass}, we show the average percentage of met demand for a different number of drones when we use/do not use demand classification.  The red color shows the percentage for the different number of drones when demand is not classified. In this model, the state of charge of batteries is either full or empty, and full batteries are used to satisfy the demand without considering the classified distance between the station and hospitals. \textcolor{black}{The system's state is the number of fully-charged batteries. The action is the number of recharging actions that make an empty battery fully-charged. The uncertainty is the stochastic demand disregarding the distance between the station and hospitals (demand nodes).} Hence, when demand is not classified, drones that return from (even a close) delivery mission are not available before recharging for the next delivery task.  For this model, optimal policies indicate that more than 150 drones are needed to satisfy 100\% of demand over 500 sample paths. In the model with demand classification, we note that finding the optimal policy and, in turn, the average percentage of demand for $M > 21$ is beyond our computational resources. Therefore, we only show the percentage derived from BI for $M \le 21$ using the color blue. However, using D$\epsilon$RL enables us to offer near-optimal policies for this model. As shown, D$\epsilon$RL (black line in Figure \ref{ClassVsNoClass}) provides an upper bound for the optimal number of drones needed to satisfy a particular level of demand for the model with demand classification. As shown, D$\epsilon$RL's policies constantly outperform the optimal policy of the model with no demand classification in terms of the average percentage of met demand. For instance, D$\epsilon$RL's policies can satisfy 100\% of the met demand with only 54 drones, which is significantly lower than the required number of batteries to hit this target when the demand is not classified. 

\section{Conclusion} \label{sec: Conc}

In this research, we addressed managing distribution operations of a drone swap station located at a drone hub to maximize the amount of stochastic met demand for flights delivering medical supplies in Rwanda, Africa. We denoted this problem as stochastic scheduling and allocation problem with multiple classes of demand (SA-MCD) where the stochastic demand is classified based on the distances a drone can fly, which is linked to the level of charge inside the drone's battery. We formulated the problem as a Markov Decision Process (MDP) model wherein the optimal policies determine the number of recharging batteries from one level to a higher level of charge over time when encountering stochastic demand from different demand classes. We solved the problem using backward induction (BI) and observed that we run into time/memory issues when the number of drones is greater than 21. Hence, we applied a reinforcement learning method with a descending $\epsilon$-greedy exploration feature (D$\epsilon$RL) to find high-quality approximate solutions quickly and overcome the time/memory issue. We designed a set of experiments to show the high performance of our D$\epsilon$RL method and obtain insights about how to manage the operations in the station to maximize the expected total weighted met demand when the model parameters vary. 

We found plenty of directions and opportunities related to this work for future research. For instance, we did not consider the difference between the length of time for different charging levels (still, our model outperformed the model with no demand classification). Future research should consider the time difference between different recharging actions to capture the system's behavior more realistically. As our transition probability function is large-scale and complex, it is worth investigating various model-free approaches, like RL/ADP/Q-learning, that can circumvent the complexities of the burdensome function, particularly when more than three classes are needed. Hence, this research, which introduces a new set of problems that can be used in many applications, opens the avenue for future studies from many perspectives. In terms of modeling and application, we can have multiple demand classification criteria, such as level of emergency (plus distance). Future research can add backlogging unsatisfied demands if it suits the application. 
Additionally, future research should consider how the operational charging and use actions for different demand classes impacts battery degradation wherein excessive charging should be avoided to lead to longer battery lifecycles.

\section*{Acknowledgement}
This research is supported by the Arkansas High Performance Computing Center which is funded through multiple National Science Foundation grants and the Arkansas Economic Development Commission.

\newpage
\bibliographystyle{unsrt}  
\bibliography{TRE}

\begin{thebibliography}{10}

\bibitem{DubaiDrone}
Ben Mutzabaugh.
\newblock Drone taxis? {D}ubai plans roll out of self-flying pods, February
  2017.
\newblock Last accessed on September 23, 2021 at
  \url{https://www.usatoday.com/story/travel/flights/todayinthesky/2017/02/13/dubai-passenger-carrying-drones-could-flying-july/97850596/}.

\bibitem{AgDrones2}
Jeremy Jensen.
\newblock Agricultural drones: How drones are revolutionizing agriculture and
  how to break into this booming market.
\newblock {UAV} Coach, April 2019.
\newblock Last accessed on September 21, 2021 at:
  \url{https://uavcoach.com/agricultural-drones/}.

\bibitem{UPSDrones}
Elizabeth Weise.
\newblock {UPS} tested launching a drone from a truck for deliveries.
\newblock USA Today, February 2017.
\newblock Last accessed on September 21, 2021 at
  \url{https://www.usatoday.com/story/tech/news/2017/02/21/ups-delivery-top-of-van-drone-workhorse/98057076/}.

\bibitem{DHL}
DHL.
\newblock Successful trial integration of {DHL} parcelcopter into logistics
  chain, May 2016.
\newblock Last accessed on September 23, 2021 at
  \url{http://www.dhl.com/en/press/releases/releases_2016/all/parcel_ecommerce/successful_trial_integration_dhl_parcelcopter_logistics_chain.html}.

\bibitem{Dhote2020}
Jeremy Dhote and Sabine Limbourg.
\newblock {Designing unmanned aerial vehicle networks for biological material
  transportation -- The case of {Brussels}}.
\newblock {\em Computers and Industrial Engineering}, 148:106652, oct 2020.

\bibitem{Malawi19}
Lauren Davitt.
\newblock Long-range drones deliver medical supplies to remote areas of
  {Malawi}.
\newblock Forbes, June 2019.
\newblock Last accessed on September 23, 2021 at
  \url{https://www.forbes.com/sites/unicefusa/2019/06/19/long-range-drones-deliver-medical-supplies-to-remote-areas-of-malawi/?sh=259e1fd36add}.

\bibitem{Unicef20}
{UNICEF Supply Division}.
\newblock How drones can be used to combat {COVID-19}.
\newblock Technical report, UNICEF, 2020.

\bibitem{McNabb20}
Miriam McNabb.
\newblock Drone delivery in the era of the pandemic: From the floor at
  commercial {UAV} expo.
\newblock drone life, September 2020.
\newblock Last accessed on September 23, 2021 at
  \url{https://dronelife.com/2020/09/16/drone-delivery-in-the-era-of-the-pandemic-from-the-floor-at-commercial-uav-expo/}.

\bibitem{Kent19}
Chloe Kent.
\newblock Moving medical supplies: enter the drone.
\newblock Medical Service Network, December 2019.
\newblock Last accessed on September 23, 2021 at
  \url{https://www.medicaldevice-network.com/features/medical-supply-drones/}.

\bibitem{Zipline20}
Kim Lyons.
\newblock Zipline and {Walmart} to launch drone deliveries of health and
  wellness products, September 2020.
\newblock Last accessed on September 23, 2021 at
  \url{https://www.theverge.com/2020/9/14/21435019/zipline-walmart-drone-deliveries-healthcare-amazon}.

\bibitem{Matternet20}
Matternet.
\newblock Matternet's {M2} drone system enabling new {U.S.} hospital delivery
  network at {Wake Forest Baptist Health}, July 2020.
\newblock Last accessed on September 23, 2021 at
  \url{https://www.suasnews.com/2020/07/matternets-m2-drone-system-enabling-new-u-s-hospital-delivery-network-at-wake-forest-baptist-health/}.

\bibitem{MannaAero20}
Simon Chandler.
\newblock Coronavirus delivers `world's first' drone delivery service.
\newblock Forbes, April 2020.
\newblock Last accessed on September 23, 2021 at
  \url{https://www.forbes.com/sites/simonchandler/2020/04/03/coronavirus-delivers-worlds-first-drone-delivery-service/?sh=407e72774957}.

\bibitem{Matternet20-2}
Mike Ball.
\newblock Matternet unveils new medical drone payload exchange station.
\newblock Unmanned Systems News, March 2020.
\newblock Last accessed on September 23, 2021 at
  \url{https://www.unmannedsystemstechnology.com/2020/03/matternet-unveils-new-medical-drone-payload-exchange-station/}.

\bibitem{Lacey13}
G.~Lacey, T.~Jiang, G.~Putrus, and R.~Kotter.
\newblock The effect of cycling on the state of health of the electric vehicle
  battery.
\newblock In {\em 48th International Universities' Power Engineering
  Conference}, pages 1--7, Dublin, Ireland, September 2013.

\bibitem{Shirk15}
Matthew Shirk and Jeffrey Wishart.
\newblock Effects of electric vehicle fast charging on battery life and vehicle
  performance.
\newblock In {\em SAE Technical Paper}, pages 1--13, Detroit, MI, April 2015.
  SAE 2015 World Congress {\&} Exhibition, SAE International.

\bibitem{electrek2018}
Fred Lambert.
\newblock {NIO} deploys 18 battery swap stations covering 2,000+ km expressway.
\newblock electrek, November 2018.
\newblock Last accessed on September 23, 2021 at
  \url{https://electrek.co/2018/11/15/nio-battery-swap-stations-network/}.

\bibitem{ChinaDaily19}
Li~Fusheng.
\newblock {BJEV} advocates battery swap service for e-vehicle owners.
\newblock China Daily, September 2019.
\newblock Last accessed on September 23, 2021 at
  \url{https://www.chinadaily.com.cn/a/201909/09/WS5d75e342a310cf3e3556a816.html}.

\bibitem{Gogoro}
{BBC}.
\newblock Electric scooter with swappable batteries hits market, June 2015.
\newblock Last accessed on September 23, 2021 at
  \url{http://www.bbc.com/news/technology-33183031}.

\bibitem{FuelRod}
FuelRod.
\newblock Power ready to go, June 2017.
\newblock Last accessed on September 23, 2021 at \url{http://www.fuel-rod.com}.

\bibitem{ChinaSwap}
Li~Dongmei.
\newblock China's {BAIC} group launches {EV} battery-swap station network in
  {Beijing}.
\newblock China Money Network, November 2016.
\newblock Last accessed on September 23, 2021 at
  \url{https://www.chinamoneynetwork.com/2016/11/03/chinas-baic-group-launches-ev-battery-swap-station-network-in-beijing}.

\bibitem{FranceSwap}
Audrey Chauvet.
\newblock S\'{e}gol\`{e}ne royal: An agreement at {COP21} would accelerate the
  energy transition, March 2015.
\newblock Last accessed on September 23, 2021 at
  \url{http://www.20minutes.fr/planete/1742767-20151203-segolene-royal-accord-cop21-permettrait-accelerer-transition-energetique}.

\bibitem{IndiaSwap}
Alnoor Peermohamed.
\newblock Reva founder looks to make electric cars affordable in {I}ndia.
\newblock Business Standard, April 2017.
\newblock Last accessed on September 23, 2021 at
  \url{http://www.business-standard.com/article/companies/reva-founder-looks-to-make-electric-cars-affordable-in-india-117041201118_1.html}.

\bibitem{SlovakiaSwap}
Nikki Gordon-Bloomfield.
\newblock Forget better place: Simple battery swap technology from {Slovakia}
  proves its worth, February 2014.
\newblock Last accessed on September 23, 2021 at
  \url{https://transportevolved.com/2014/02/12/forget-better-place-simple-battery-swap-technology-from-slovakia-proves-its-worth/}.

\bibitem{Puterman05}
Martin~L. Puterman.
\newblock {\em Markov decision processes: Discrete stochastic dynamic
  programming}.
\newblock John Wiley \& Sons, Hoboken, New Jersey, 1st edition, 2005.

\bibitem{Powell}
Warren~B. Powell.
\newblock {\em Approximate Dynamic Programming: Solving the Curses of
  Dimensionality}.
\newblock John Wiley \& Sons, New York, NY, USA, 2 edition, 2011.

\bibitem{ZiplineCarry}
Aryn Baker.
\newblock The {American} drones saving lives in {Rwanda}, July 2017.
\newblock Last accessed on September 23, 2021 at
  \url{https://time.com/rwanda-drones-zipline/}.

\bibitem{Asadi19}
Amin Asadi and Sarah {Nurre Pinkley}.
\newblock A stochastic scheduling, allocation, and inventory replenishment
  problem for battery swap stations.
\newblock {\em Transportation Research Part E: Logistics and Transportation
  Review}, 146:102212, 2021.

\bibitem{AsadiTs2021}
Amin Asadi and Sarah {Nurre Pinkley}.
\newblock A monotone approximate dynamic programming approach for the
  stochastic scheduling, allocation, and inventory replenishment problem:
  Applications to drone and electric vehicle battery swap stations.
\newblock {\em CoRR}, abs/2105.07026, 2021.

\bibitem{Widrick16}
Rebecca~S. Widrick, Sarah~G. Nurre, and Matthew~J. Robbins.
\newblock Optimal policies for the management of an electric vehicle battery
  swap station.
\newblock {\em Transportation Science}, 52(1):59--79, 2018.

\bibitem{Nurre14}
Sarah~G. Nurre, Russell Bent, Feng Pan, and Thomas~C. Sharkey.
\newblock Managing operations of plug-in hybrid electric vehicle {(PHEV)}
  exchange stations for use with a smart grid.
\newblock {\em Energy Policy}, 67:364--377, 2014.

\bibitem{Kwizera2018}
Olivier Kwizera and Sarah~G. Nurre.
\newblock Using drones for delivery: A two-level integrated inventory problem
  with battery degradation and swap stations.
\newblock In {\em Proceedings of the Industrial and Systems Engineering
  Research Conferences}, pages 1--6, Orlando, FL, 2018.

\bibitem{Macrina20}
Giusy Macrina, Luigi {Di Puglia Pugliese}, Francesca Guerriero, and Gilbert
  Laporte.
\newblock {Drone-aided routing: A literature review}.
\newblock {\em Transportation Research Part C: Emerging Technologies},
  120:102762, nov 2020.

\bibitem{Barmpounakis16}
Emmanouil~N. Barmpounakis, Eleni~I. Vlahogianni, and John~C. Golias.
\newblock {Unmanned Aerial Aircraft Systems for transportation engineering:
  Current practice and future challenges}.
\newblock {\em International Journal of Transportation Science and Technology},
  5(3):111--122, oct 2016.

\bibitem{Chang2018}
Yong~Sik Chang and Hyun~Jung Lee.
\newblock {Optimal delivery routing with wider drone-delivery areas along a
  shorter truck-route}.
\newblock {\em Expert Systems with Applications}, 104:307--317, aug 2018.

\bibitem{Otto2018}
Alena Otto, Niels Agatz, James Campbell, Bruce Golden, and Erwin Pesch.
\newblock {Optimization approaches for civil applications of unmanned aerial
  vehicles (UAVs) or aerial drones: A survey}.
\newblock {\em Networks}, 72(4):411--458, 2018.

\bibitem{Khoufi2019}
Ines Khoufi, Anis Laouiti, and Cedric Adjih.
\newblock {A survey of recent extended variants of the traveling salesman and
  vehicle routing problems for unmanned aerial vehicles}.
\newblock {\em Drones}, 3(3):1--30, 2019.

\bibitem{Seeker}
Tracy Staedter.
\newblock Drones now delivering life-saving blood in {Rwanda}, October 2016.
\newblock Last accessed on September 21, 2021 at
  \url{https://www.seeker.com/drones-deliver-blood-rwanda-2045341414.html}.

\bibitem{Rescue2019}
Robyn Bainbridge.
\newblock Special report: {Zipline} international blood drone delivery service,
  August 2019.
\newblock Last accessed on September 21, 2021 at
  \url{https://www.airmedandrescue.com/latest/long-read/special-report-zipline-international-blood-drone-delivery-service}.

\bibitem{Vincent21}
James Vincent.
\newblock {Self-flying drones are helping speed deliveries of COVID-19 vaccines
  in Ghana}, May 2021.
\newblock Last accessed on September 21, 2021 at
  \url{https://www.theverge.com/2021/3/9/22320965/drone-delivery-vaccine-ghana-zipline-cold-chain-storage}.

\bibitem{Draganfly}
Ishveena Singh.
\newblock {Draganfly to start drone delivery of COVID-19 vaccine to rural
  Texas}, May 2021.
\newblock Last accessed on September 21, 2021 at
  \url{https://dronedj.com/2021/05/20/draganfly-drone-delivery-covid-19/}.

\bibitem{Tokekar13}
Pratap Tokekar, Joshua Vander~Hook, David Mulla, and Volkan Isler.
\newblock Sensor planning for a symbiotic {UAV} and {UGV} system for precision
  agriculture.
\newblock In {\em 2013 IEEE/RSJ International Conference on Intelligent Robots
  and Systems}, pages 5321--5326, 2013.

\bibitem{Wang2017}
Xingyin Wang, Stefan Poikonen, and Bruce Golden.
\newblock The vehicle routing problem with drones: several worst-case results.
\newblock {\em Optimization Letters}, 11(4):679--697, 2017.

\bibitem{Savuran2016}
Halil Savuran and Murat Karakaya.
\newblock Efficient route planning for an unmanned air vehicle deployed on a
  moving carrier.
\newblock {\em Soft Computing}, 20(7):2905--2920, 2016.

\bibitem{Guerriero2014}
F.~Guerriero, R.~Surace, V.~Loscr{\'\i}, and E.~Natalizio.
\newblock A multi-objective approach for unmanned aerial vehicle routing
  problem with soft time windows constraints.
\newblock {\em Applied Mathematical Modelling}, 38:839--852, 2 2014.

\bibitem{ZiplineFeatures}
{Engineering for Change}.
\newblock Special report: {Zipline} international blood drone delivery service,
  May 2021.
\newblock Last accessed on September 21, 2021 at
  \url{https://www.engineeringforchange.org/solutions/product/zipline/}.

\bibitem{Al-Sabban13}
W.~H. {Al-Sabban}, L.~F. {Gonzalez}, and R.~N. {Smith}.
\newblock Wind-energy based path planning for {Unmanned Aerial Vehicles} using
  {Markov Decision Processes}.
\newblock In {\em 2013 IEEE International Conference on Robotics and
  Automation}, pages 784--789, May 2013.

\bibitem{Baek13}
S.~S. {Baek}, H.~{Kwon}, J.~A. {Yoder}, and D.~{Pack}.
\newblock Optimal path planning of a target-following fixed-wing {UAV} using
  sequential decision processes.
\newblock In {\em 2013 IEEE/RSJ International Conference on Intelligent Robots
  and Systems}, pages 2955--2962, November 2013.

\bibitem{Yu15}
Y.~{Fu}, X.~{Yu}, and Y.~{Zhang}.
\newblock Sense and collision avoidance of unmanned aerial vehicles using
  markov decision process and flatness approach.
\newblock In {\em Proceeding of the 2015 IEEE International Conference on
  Information and Automation}, pages 714--719, August 2015.

\bibitem{Federgruen84}
Awi Federgruen and Paul Zipkin.
\newblock Approximations of dynamic, multilocation production and inventory
  problems.
\newblock {\em Management Science}, 30(1):69--84, 1984.

\bibitem{SOMARIN17}
Aghil~Rezaei Somarin, Songlin Chen, Sobhan Asian, and David~Z.W. Wang.
\newblock A heuristic stock allocation rule for repairable service parts.
\newblock {\em International Journal of Production Economics}, 184:131--140,
  2017.

\bibitem{millart04}
Oguzhan Alagoz, Lisa~M. Maillart, Andrew~J. Schaefer, and Mark~S. Roberts.
\newblock The optimal timing of living-donor liver transplantation.
\newblock {\em Management Science}, 50(10):1420--1430, 2004.

\bibitem{Chhatwal10}
Jagpreet Chhatwal, Oguzhan Alagoz, and Elizabeth~S. Burnside.
\newblock Optimal breast biopsy decision-making based on mammographic features
  and demographic factors.
\newblock {\em Operations Research}, 58(6):1577--1591, 2010.

\bibitem{Zhang12}
Jingyu Zhang, Brian~T. Denton, Hari Balasubramanian, Nilay~D. Shah, and
  Brant~A. Inman.
\newblock Optimization of prostate biopsy referral decisions.
\newblock {\em Manufacturing \& Service Operations Management}, 14(4):529--547,
  2012.

\bibitem{Khojandi14}
Anahita Khojandi, Lisa~M. Maillart, Oleg~A. Prokopyev, Mark~S. Roberts, Timothy
  Brown, and William~W. Barrington.
\newblock Optimal implantable cardioverter defibrillator {(ICD)} generator
  replacement.
\newblock {\em INFORMS Journal on Computing}, 26(3):599--615, 2014.

\bibitem{Gayon09}
Jean-Philippe Gayon, Saif Benjaafar, and Francis de~V{\'e}ricourt.
\newblock Using imperfect advance demand information in production-inventory
  systems with multiple customer classes.
\newblock {\em Manufacturing \& Service Operations Management}, 11(1):128--143,
  2009.

\bibitem{Benjaafar2011}
Saif Benjaafar, Mohsen ElHafsi, Chung~Yee Lee, and Weihua Zhou.
\newblock {Optimal control of an assembly system with multiple stages and
  multiple demand classes}.
\newblock {\em Operations Research}, 59(2):522--529, 2011.

\bibitem{Thompson2009}
Steven Thompson, Manuel Nunez, Robert Garfinkel, and Matthew~D. Dean.
\newblock {Efficient short-term allocation and reallocation of patients to
  floors of a hospital during demand surges}.
\newblock {\em Operations Research}, 57(2):261--273, 2009.

\bibitem{Milnar16}
Tanja Mlinar and Philippe Chevalier.
\newblock {Dynamic admission control for two customer classes with stochastic
  demands and strict due dates}.
\newblock {\em International Journal of Production Research},
  54(20):6156--6173, 2016.

\bibitem{WANG2019}
Hua Wang, De~Zhao, Qiang Meng, Ghim~Ping Ong, and Der-Horng Lee.
\newblock A four-step method for electric-vehicle charging facility deployment
  in a dense city: An empirical study in {Singapore}.
\newblock {\em Transportation Research Part A: Policy and Practice},
  119:224--237, 2019.

\bibitem{Sutton18}
Richard~S. Sutton and Andrew~G. Barto.
\newblock {\em Reinforcement Learning: An Introduction}.
\newblock The MIT Press, Cambridge, MA, 2 edition, 2018.

\bibitem{vanr97}
B.~Van Roy, D.P. Bertsekas, Y.~Lee, and J.N. Tsitsiklis.
\newblock A neuro-dynamic programming approach to retailer inventory
  management.
\newblock In {\em Proceedings of the {IEEE} Conference on Decision and
  Control}, volume~4, pages 4052--4057, 1997.

\bibitem{Cimen17}
Mustafa {\c C}imen and Chris Kirkbride.
\newblock Approximate dynamic programming algorithms for multidimensional
  flexible production-inventory problems.
\newblock {\em International Journal of Production Research}, 55(7):2034--2050,
  2017.

\bibitem{Cimen15}
Mustafa {\c C}imen and Christopher Kirkbride.
\newblock Approximate dynamic programming algorithms for multidimensional
  inventory optimization problems.
\newblock {\em IFAC Proceedings Volumes}, 46(9):2015--2020, 2013.

\bibitem{JiangSheng09}
Chengzhi Jiang and Zhaohan Sheng.
\newblock Case-based reinforcement learning for dynamic inventory control in a
  multi-agent supply-chain system.
\newblock {\em Expert Systems with Applications}, 36:6520--6526, 2009.

\bibitem{Chaharsooghi}
S.~Kamal Chaharsooghi, Jafar Heydari, and S.~Hessameddin Zegordi.
\newblock A reinforcement learning model for supply chain ordering management:
  An application to the beer game.
\newblock {\em Decision Support Systems}, 45(4):949--959, 2008.

\bibitem{Bertmis02}
Dimitris Bertsimas and Ramazan Demir.
\newblock An approximate dynamic programming approach to multidimensional
  knapsack problems.
\newblock {\em Management Science}, 48(4):550--565, 2002.

\bibitem{Erdelyi10}
Alexander Erdelyi and Huseyin Topaloglu.
\newblock Approximate dynamic programming for dynamic capacity allocation with
  multiple priority levels.
\newblock {\em IIE Transactions}, 43(2):129--142, 2010.

\bibitem{Maxwell10}
Matthew~S. Maxwell, Mateo Restrepo, Shane~G. Henderson, and Huseyin Topaloglu.
\newblock Approximate dynamic programming for ambulance redeployment.
\newblock {\em INFORMS Journal on Computing}, 22(2):266--281, 2010.

\bibitem{Powell05}
Warren~B. Powell and Huseyin Topaloglu.
\newblock Approximate dynamic programming for large-scale resource allocation
  problems.
\newblock {\em INFORMS Tut{OR}ials in Operations Research}, pages 123--147,
  2005.

\bibitem{Nasrollahzadeh18}
Amirali Nasrollahzadeh, Amin Khademi, and Maria~E. Mayorga.
\newblock Real-time ambulance dispatching and relocation.
\newblock {\em Manufacturing and Service Operations Management},
  20(3):467--480, 2018.

\bibitem{Kwon08}
Ick-Hyun Kwon, Chang~Ouk Kim, Jin Jun, and Jung~Hoon Lee.
\newblock Case-based myopic reinforcement learning for satisfying target
  service level in supply chain.
\newblock {\em Expert Systems with Applications}, 35(1):389--397, 2008.

\bibitem{Ryzhov19}
Ilya~O. Ryzhov, Martijn R.~K. Mes, Warren~B. Powell, and Gerald van~den Berg.
\newblock {Bayesian exploration strategies for approximate dynamic
  programming}.
\newblock {\em Operations Research}, 67(1):198--214, 2019.

\bibitem{Zipline80km}
Evan Ackerman.
\newblock Zipline launches long-distance drone delivery of {COVID-19} supplies
  in the {U.S.}, May 2020.
\newblock Last accessed on September 21, 2021 at
  \url{https://spectrum.ieee.org/automaton/robotics/drones/zipline-long-distance-delivery-covid19-supplies}.

\bibitem{Sinnott84}
R.W. Sinnott.
\newblock Virtues of the {Haversine}.
\newblock {\em Sky and Telescope}, 68(2):158, 1984.

\bibitem{DroneAirport}
{Rwanda Civil Aviation Authority}.
\newblock Unmanned aircraft operations in {Rwanda} unmanned aircraft operations
  in {Rwanda}, May 2021.
\newblock Last accessed on September 21, 2021 at
  \url{https://caa.gov.rw/index.php?id=110}.

\bibitem{Swartzman70}
Gordon Swartzman.
\newblock The patient arrival process in hospitals: statistical analysis.
\newblock {\em Health services research}, 5(4):320--329, Winter 1970.

\bibitem{Armony2015}
Mor Armony, Shlomo Israelit, Avishai Mandelbaum, Yariv~N. Marmor, Yulia
  Tseytlin, and Galit~B. Yom-Tov.
\newblock {On patient flow in hospitals: A data-based queueing-science
  perspective}.
\newblock {\em Stochastic Systems}, 5(1):146--194, 2015.

\bibitem{DroneRecharge}
Harry {McNabb}.
\newblock A drone battery that charges in 5 minutes, August 2020.
\newblock Last accessed on September 21, 2021 at
  \url{https://dronelife.com/2020/08/01/a-drone-battery-that-charges-in-5-minutes/}.

\bibitem{DroneSpeed}
Magdalena Petrova and Lora Kolodny.
\newblock Zipline's new drone can deliver medical supplies at 79 miles per
  hour, April 2018.
\newblock Last accessed on September 21, 2021 at
  \url{https://www.cnbc.com/2018/04/02/zipline-new-zip-2-drone-delivers-supplies-at-79-mph.html}.

\bibitem{BloodEstimate}
Neelam Dhingra.
\newblock Estimate blood requirements - search for a global standard, April
  2010.
\newblock Last accessed on September 21, 2021 at
  \url{https://www.who.int/bloodsafety/transfusion_services/estimation_present}.

\bibitem{BloodOverEstimate}
Ange Iliza.
\newblock How can {Rwanda} tackle its blood donation shortfall?, June 2020.
\newblock Last accessed on September 21, 2021 at
  \url{https://www.newtimes.co.rw/news/rwanda-tackle-its-blood-donation-shortfall}.

\bibitem{Green07}
Linda Green, Peter Kolesar, and Ward Whitt.
\newblock Coping with time-varying demand when setting staffing requirements
  for a service system.
\newblock {\em Production and Operations Management}, 16:13--39, 01 2007.

\bibitem{Tiwari2014}
Yogesh Tiwari, Sonu Goel, and Amarjeet Singh.
\newblock {Arrival time pattern and waiting time distribution of patients in
  the emergency outpatient department of a tertiary level health care
  institution of North India}.
\newblock {\em Journal of Emergencies, Trauma and Shock}, 7(3):160--165, 2014.

\bibitem{Jones2007}
Spencer~S Jones, M~Stat, and Todd Allen.
\newblock {Department Cycle to Improve Patient Flow}.
\newblock {\em Radiology}, 33(5):247--255, 2007.

\bibitem{George06}
Abraham~P. George and Warren~B. Powell.
\newblock Adaptive stepsizes for recursive estimation with applications in
  approximate dynamic programming.
\newblock {\em Machine Learning}, 65(1):167--198, 2006.

\end{thebibliography}

\end{document}